\documentclass{amsart}[14pt]

\usepackage{bm}
\usepackage{amsmath, amstext, amssymb, amsopn, amsthm, amsgen, amscd, 
 fouriernc, upgreek}
\usepackage{stmaryrd}
\usepackage[small,nohug,heads=vee]{diagrams} \diagramstyle[labelstyle=\scriptstyle]

\usepackage{color}
\usepackage{tikz}
\usepackage{float}

\usepackage{graphicx}
\usepackage{epstopdf}

\definecolor{CadetBlue}{cmyk}{0.62, 0.57, 0.23, 0 }
\definecolor{black}{cmyk}{1, 0.5, 0, 0 }
\definecolor{RedViolet}{cmyk}{0.07, 0.9, 0, 0.34 }
\definecolor{SeaGreen}{cmyk}{0.69, 0, 0.5, 0}

\DeclareMathAlphabet{\mathpzc}{OT1}{pzc}{m}{it}

\newcommand{\R}{\mathbb R}
\newcommand{\C}{\mathbb C}

\newcommand{\F}{\mathbb F}
\newcommand{\HP}{\mathbb H}

\newcommand{\N}{\mathbb N}
\newcommand{\Q}{\mathbb Q}

\newcommand{\Z}{\mathbb Z}

\newcommand{\T}{\mathbb T}

\newcommand{\SI}{\mathbb S}

\newtheorem{theo}{Theorem}
\newtheorem{lemm}{Lemma}
\newtheorem{prop}{Proposition}
\newtheorem{coro}{Corollary}

\newtheorem*{conj}{Conjecture}

\theoremstyle{definition}
\newtheorem{defi}{Definition}

\theoremstyle{remark}

\newtheorem{note}{Note}

\newcommand{\bast}{{}^{\ast}}

\newcommand{\bdiam}{{}^{\diamond}}
\newcommand{\bbull}{{}^{\bullet}}

\newcommand{\kf}{\mathcal{F}(\upmu ,\uptheta)}

\title{Modular Invariant of Quantum Tori}
\author{C. Casta\~{n}o Bernard \& T. M. Gendron}
\address{Centro de Estudios en F\'{i}sica y Matem\'{a}ticas B\'{a}sicas y Aplicadas, Universidad Aut\'{o}noma de Chiapas,
4a. Oriente Norte No. 1428, Colonia Barrio la Pimienta, Tuxtla Guti\'{e}rrez, Chiapas, M\'{e}xico }
\email{ccastanobernard@gmail.com}
\address{Instituto de Matem\'{a}ticas -- Unidad Cuernavaca, Universidad
Nacional Autonoma de M\'{e}xico, Av. Universidad S/N, C.P. 62210
Cuernavaca, Morelos, M\'{e}xico}
\email{tim@matcuer.unam.mx}
\date{31 January 2012}
\begin{document}
\vspace{2cm} \maketitle

\begin{abstract}  The quantum modular invariant $j^{\rm qt}(\uptheta )$ of $\uptheta\in\R$ 
is defined as a discontinuous ${\rm PGL}_{2}(\Z )$-invariant multi-valued map using the distance-to-the-nearest-integer function $\|\cdot \|$.   For $\uptheta\in\Q$ it is shown that $j^{\rm qt}(\uptheta )=\infty$ and for quadratic irrationalities PARI/GP experiments suggest that $j^{\rm qt}(\uptheta )$ is a finite set.
In the case of the golden
mean $\upvarphi$, we produce explicit formulas involving weighted
versions of the Rogers-Ramanujan functions for the experimental supremum and 
 infimum. 
 We then define 
a universal modular invariant $\bdiam j: \bdiam \widehat{\sf Mod}\rightarrow \bdiam \hat{\C}$
as a {\it continuous} and {\it single valued} map of ultrasolenoids, such that 1) the classical modular invariant 
is a quotient of the restriction of $\bdiam j$ to a subsolenoid  ${\sf Mod}^{\rm cl}\subset \bdiam \widehat{\sf Mod}$
fibering over the classical moduli space of elliptic curves
and 2) the quantum modular invariant
is a quotient of the restriction of $\bdiam j$ to a subsolenoid  ${\sf Mod}^{\rm qt}\subset \bdiam \widehat{\sf Mod}$
fibering over the moduli space of elliptic curves equipped with a Kronecker foliation.

\end{abstract}
\tableofcontents

\section*{Introduction}

One of the most powerful links between between number theory, geometry and analysis can be found in
the theory of complex multiplication (CM): a signal event in a long development in number theory,
which begins with Gau\ss's reciprocity law and
culminates in the main theorems of class field theory \cite{Se}, \cite{Si}.

The theory of CM yields
results that come under the heading of \guillemotleft explicit class field theory\guillemotright, which may be seen as generalizations of the Theorem of Kronecker-Weber.  If $\upomega\in\C-\Q$ is a complex quadratic irrationality, let $E_{\upomega}(\C )\cong\C/\langle1,\upomega\rangle $ be the complex points of the
elliptic curve $E_{\upomega}$ parametrized by $\upomega$.  Then $E_{\upomega}(\C )$ has CM (its endomorphism ring is strictly larger than $\Z$); if we denote $K=\Q (\upomega )$ and assume that ${\rm End}(E_{\upomega}(\C))\otimes\Q =K$ then
\begin{itemize}
\item The maximal unramified abelian extension (Hilbert class field) $H$ of $K$ is generated over $K$ by the modular invariant $j(\upomega )$ of 
$E_{\upomega}(\C )$.
\item The maximal abelian extension $K^{\rm ab}$ of $K$ is generated over $H$ by the values of the Weierstra\ss\ $\wp$-function
at the torsion subgroup of $\C/\langle 1,\upomega\rangle$.
\end{itemize}

Finding the analogue of the theory of CM for more general fields -- Kronecker's {\it Jugendtraum} or Hilbert's twelfth problem -- has 
been the focus of investigation for nearly a century \cite{He}, \cite{Sh}, \cite{St}, \cite{Ka}, \cite{DDP}.  

In 2004, Yu.\ Manin
\cite{Ma} proposed an approach to the Stark conjectures \cite{St}
in which
quantum tori play a role analogous to
that of elliptic curves with CM.
Manin's {\it Alterstraum} is sometimes known
as the \guillemotleft Real Multiplication\guillemotright\ programme (RM): an approach to the Stark conjectures in the case of real quadratic fields which uses notions of noncommutative geometry.

In \cite{Ma}, the quantum torus $\T(\uptheta)$,  $\uptheta\in\R -\Q$, is understood as an object in a category
of irrational rotation $C^{\ast}$-algebras, however
it may also be described in somewhat more naive geometric terms as 
the quotient $\R/\langle 1,\uptheta \rangle$ of the reals by the
pseudo lattice $\langle 1,\uptheta \rangle$, or equivalently, as
the space of leaves of 
the Kronecker foliation $\mathcal{F}(\uptheta)$ {\rm c.f.}  \cite{Mar}
or \S \ref{qtoriandkfol} of this paper.  The moduli space of quantum tori is identified with the quotient 
${\sf Mod}^{\rm qt}={\rm PGL}_{2}(\Z )\backslash (\R-\Q)$.

One would therefore like to formulate and prove exact analogues
of the main theorems of CM in the RM case, 
using suitable notions of  Weierstra\ss\ $\wp$-function and modular invariant for the quantum torus $\T(\uptheta )$.
However both $\T(\uptheta )$ and ${\sf Mod}^{\rm qt}$ are quotients by groups acting with dense orbits, so
in particular, there are no non constant continuous functions defined on either of them.
Thus it is not at all clear how to define the analogues
of the Weierstra\ss\ $\wp$-function or modular invariant in this setting.


The goal of the present paper is to provide a definition of the modular invariant of quantum tori:  more precisely we 
\begin{enumerate}
\item[A.] Give an elementary definition of the (necessarily discontinuous) quantum modular invariant $j^{\rm qt}(\uptheta)$ of a real number $\uptheta\in\R$, using
only the distance-to-the-nearest-integer function.   The association $\uptheta\mapsto j^{\rm qt}(\uptheta)$ induces a 
 {\it multi-valued} function of ${\sf Mod}^{\rm qt}$ i.e. taking values in ${\sf 2}^{\R}\cup\{\infty\}$,
 which could be interpreted as the spectrum of an (as yet to be discovered) operator.  Experimental evidence suggests that the set $j^{\rm qt}(\uptheta)$ is finite if $\uptheta$ is a quadratic irrationality.

\item[B.] Define a {\it universal} modular invariant 
\[ \bdiam j :  \bdiam \widehat{\sf Mod}\rightarrow \bdiam \hat{\C}\]
as a {\it continuous} and {\it single valued} map of {\it ultrasolenoids} (see further below or \S \ref{ultrasolenoids}) and show that both the classical modular invariant $j^{\rm cl}$ and the quantum modular invariant $j^{\rm qt}$
occur are subquotients of $\bdiam j$.  The ultrasolenoid $\bdiam \widehat{\sf Mod}$ can thus be construed as the ``Riemann surface'' associated to the
multi-valued and discontinuous $j^{\rm qt}$. 
\end{enumerate}

We describe in more detail how the above is accomplished by way of an overview of the sections making up this paper.

In \S \ref{quantummodinvar} we give the definition of $j^{\rm qt}(\uptheta )\subset \R\cup \{\infty\}$, show that it is invariant 
with respect to the action of ${\rm PGL}_{2}(\Z )$ 
and that for $\uptheta\in\Q$, $j^{\rm qt}(\uptheta )=\infty$.
 In \S\S \ref{golden},\ref{golden2} we discuss the case of $\uptheta = \upvarphi$ = the golden mean.
We deduce an {\it explicit formula} for a value $ j^{\rm qt} (\upvarphi)$ which experiment suggests is the infimum of $ j^{\rm qt} (\upvarphi)$
(as well as a slight modification which appears to coincide with the experimental supremum of $j^{\rm qt}(\upvarphi )$) in terms of weighted variants of the Rogers-Ramanujan functions. This concludes the part of the paper corresponding to A. above, the  \guillemotleft elementary\guillemotright\ part of the paper.

The remaining sections are devoted to removing the discontinuity and multivaluedness of $j^{\rm qt}$ by extending and lifting
it to a larger space modeled on the Anosov foliation.  The way that this is accomplished has the added benefit 
of putting the invariant into a universal geometric context, allowing us to argue in favor of its interpretation as
the modular invariant of quantum tori, as well as relate it to the classical modular invariant.  

In \S 4 
we define the {\it generalized Kronecker foliation} $\mathcal{F}(\upmu,\uptheta )$ of the elliptic curve $\T(\upmu)$ by lines of
$\upmu$-slope $\uptheta$ and establish its relation to the quantum torus $\T(\uptheta )$.  We prove that the moduli space
 of generalized Kronecker foliations is the  \guillemotleft Anosov foliation\guillemotright\
\[ {\sf Mod}^{\rm kf} \approx {\rm PGL}_{2}(\Z)\backslash(\pm\HP\times( \R\cup\{\infty\})),\]
where $\pm\HP=\HP\cup\overline{\HP}$.  The leaf space of ${\sf Mod}^{\rm kf}$ is the completed moduli space of quantum tori
$\overline{\sf Mod}^{\rm qt}= {\sf Mod}^{\rm qt}\cup\{\infty\}$.

In \S \ref{nonstdstruct} we review the notions of ultrafilter and ultrapower, then define (as ultrapowers over $\N$) the nonstandard integers, reals and complexes $\bast\Z\subset \bast\R\subset\bast\C$. In
\S \ref{DAGs} we define a \guillemotleft uniformizing lattice\guillemotright\ for the Kronecker foliation $\mathcal{F}(\upmu,\uptheta )$:
the {\it diophantine approximation group} \[ \bast\Uplambda(\upmu ,\uptheta ) 
\subset\bast\Uplambda(\upmu ) \]
where $\Uplambda(\upmu )=\langle 1,\upmu\rangle$ and $\bast\Uplambda(\upmu )$ is its ultrapower.  If $\bbull\R\subset\bbull\C$ denote 
the vector spaces $\bast\R\subset\bast\C$ modulo infinitesimals then $ \bast\Uplambda(\upmu ,\uptheta ) \subset \bbull\C$ stabilizes the line $\bbull\R\cdot (1+\upmu\uptheta )$
with quotient isomorphic to $\mathcal{F}(\upmu,\uptheta )$.

The path is then clear: define Eisenstein series following the usual prescription -- but using $\bast\Uplambda(\upmu ,\uptheta )$ in place
of $\Uplambda(\upmu )$ -- to arrive at a modular invariant that will be a (necessarily {\it transversally} discontinuous) function 
of ${\sf Mod}^{\rm kf}$.  The challenge is to make sense of summation over the uncountable group $\bast\Uplambda(\upmu ,\uptheta )$, which we do by passing to a sheaf of ultrapowers of $\bast\C$.  The net effect will be an expansion of the domain and codomain
of our invariant.

In fact, it makes sense to carry out this task in a universal way which shepherds both the classical and quantum invariant 
into the confines of a single invariant.  Given $\upmu\in\pm\HP$ and $[F_{\upalpha}]\subset\bast\Z^{2}-\{0,0\}$ a {\it hyperfinite subset}, we may
form the {\it hyperfinite partial sum}
 (see
\S\S \ref{ultrasolenoids},\ref{eisection})
\[ G_{k}(\upmu )_{[F_{i}]}   = \sum_{(\bast m,\bast n)\in [F_{i}]}(\bast m\upmu+\bast n)^{-2k}\in\bast\C .\] 
One obtains a net of partial sums indexed by the set $\mathcal{H}$ of hyperfinite subsets of $\bast\Z^{2} -\{0,0\}$.  This net defines a section $\bdiam\breve{G}_{k}(\upmu )$ of the sheaf $\bdiam\breve{\C}$ of ultrapowers of $\bast\C$ over the Stone space ${\sf Ult}(\mathcal{H})$ of ultrafilters on $\mathcal{H}$.  That is, evaluating the usual formula for the modular invariant at the above-defined Eisenstein sections
and at the appropriate values of $k$, we obtain a function 
\[ \bdiam \breve{\jmath}:\pm\HP\rightarrow \bdiam\breve{\Upgamma} = \text{the set of sections of } \bdiam\breve{\C}.\]  

The function $\bdiam \breve{\jmath}$
is however {\it not} ${\rm GL}_{2}(\Z )$-invariant, since there is an attendant shift of indices when acting by elements of ${\rm GL}_{2}(\Z )$.  We can nevertheless achieve modularity   
by quotienting $\bdiam\breve{\C}$ by
the shift action of ${\rm GL}_{2}(\Z )$ on $\bdiam\breve{\C}$, producing the {\it ultrasolenoid} $\bdiam\hat{\C}$.   Then if we denote
by $\bdiam\hat{\Upgamma}={\rm GL}_{2}(\Z )\backslash\bdiam\breve{\Upgamma}=$  the $\C$-algebra of {\it ultratransversals},
we obtain a ${\rm GL}_{2}(\Z )$ invariant function
\[  \bdiam\hat{\jmath}: \pm\HP\longrightarrow \bdiam\hat{\Upgamma}.\]
This is discussed in the first part of \S \ref{modinvariantsection}.

More insightfully, if we let $\bdiam {\sf Mod}$ be the solenoid obtained as the quotient of the product  $\pm \HP\times {\sf Ult}(\mathcal{H})\subset \bdiam\breve{\C}$
by the diagonal action of ${\rm GL}_{2}(\Z)$ (shift on the base, linear action along the stalks), then we obtain, equivalently,
a universal leaf preserving, transversally {\it continuous} function
\[ \bdiam j: \bdiam {\sf Mod}\longrightarrow \bdiam\hat{\C}. \]
The space  $\bdiam {\sf Mod}$ is an obvious generalization of the Anosov foliation ${\sf Mod}^{\rm kf}$, where $\R\cup \{\infty\}$ has
been exploded and retopologized to the locally Cantor ${\sf Ult}(\mathcal{H})$.
See Theorem \ref{solvaluedinvariant} at the end of
\S \ref{modinvariantsection}.  

To recover the classical and quantum invariants, we select out ultrafilters that  \guillemotleft observe\guillemotright\ the groups
$\bast\Uplambda (\upmu)$ resp.\ $\bast\Uplambda (\upmu,\uptheta)$.  These are the ${\rm GL}_{2}(\Z )$-invariant
subspaces 
\[ {\sf Cone}^{\rm cl},\;\;   {\sf Cone}^{\rm qt} = \bigsqcup_{\uptheta\in\R} {\sf Cone}^{\rm qt}(\uptheta ) \subset {\sf Ult}(\mathcal{H})\]
of classical and quantum {\it cone ultrafilters}, see \S \ref{ultrasolenoids}.  They give rise to subultrasolenoids 
\[ \bdiam {\sf Mod}^{\rm cl}, 
\bdiam {\sf Mod}^{\rm qt}\subset \bdiam {\sf Mod},\] and the restriction of $\bdiam j$ to each defines the classical resp.\ quantum invariants $\bdiam j^{\rm cl}$, $\bdiam j^{\rm qt}$.
If we denote by $\simeq$ the relation of infinitesimality, then
\[  \bdiam j^{\rm cl} (\upmu, \mathfrak{u})\simeq j^{\rm cl}(\upmu )\quad \forall\upmu\in\pm\HP \]
moreover, every limit point $a\in j^{\rm qt}(\uptheta )$ is near standard to $ \bdiam j^{\rm qt} (i, \mathfrak{u})$
for some $ \mathfrak{u}\in  {\sf Cone}^{\rm qt}(\uptheta )$.  
These statements are proved (at the level of ultratransversal valued invariants) in Corollary \ref{classasympclass} and Theorem \ref{stdprtlim}
of \S \ref{modinvariantsection}; their rendering into the language of ultrasolenoid valued invariants is made using Theorem \ref{solvaluedinvariant}.

In the Appendix we have presented some PARI/GP calculations which suggest that at the quadratic irrationalities, $j(\uptheta )$ is a finite set.

Ours is not the first attempt to use nonstandard constructions in the consideration of the RM problem: 
the reader may wish to compare the ideas in this paper with the work of Fesenko \cite{Fe1}, \cite{Fe2} and his student Taylor
\cite{Ta1}, \cite{Ta2}, \cite{Ta3}.   Approaches using noncommutative geometry are discussed in the review \cite{Ma1}.

\vspace{5mm} 

\noindent {\bf {\small Acknowledgments:}}  We would like to thank B. Zil'ber for having
suggested that our ideas might be applicable to the problem of RM and  
R. Kossak for having brought to our attention \cite{Mac}.  In addition, we are indebted to D. Zagier for pointing out
the multivalued nature of $j^{\rm qt}$ and making useful suggestions for its experimental analysis.
The second author would like to express his gratitude to the Isaac Newton Institute for making possible a visit in April of 2013: without which his conversations with D. Zagier would not have taken place.
This work was supported in part by the grants CONACyT 058537 and PAPIIT IN103708.

\section{The Quantum Modular Invariant of a Real Number}\label{quantummodinvar}

Fix $\uptheta\in\R$.   Let $\|\cdot \|:\R\rightarrow[0,1/2]$ denote the function which assigns to a 
real number its distance to the nearest integer.  
If $n\in\Z$ satisfies $\|n\uptheta\|<1/2$ we denote by $n^{\perp}$ the unique closest integer, so that if we write
\[ \upvarepsilon (n):= n\uptheta -n^{\perp}\] then $|\upvarepsilon (n)|=\|n\uptheta\|$.

   For $\upvarepsilon >0$ let
 \[ B_{\upvarepsilon}(\uptheta ) = \big\{   n\in \N \big|\; \|n\uptheta\|<\upvarepsilon   \big\} \]
and define the {\bf {\small $\boldsymbol\upvarepsilon$ zeta function}} of $\uptheta$ as 
\[ \upzeta_{\uptheta,\upvarepsilon}(s) := \sum_{B_{\upvarepsilon}(\uptheta )}n^{-s}.\]
Define 
\[J_{\upvarepsilon}(\uptheta ) :=  \frac{49}{40} \frac{ \upzeta_{\uptheta,\upvarepsilon}(6)^{2}}{ \upzeta_{\uptheta,\upvarepsilon}(4)^{3}}
\]
and
 \begin{align}\label{stdformulaepsilon} j_{\upvarepsilon}(\uptheta ) := \frac{12^{3}}{ 1-J_{\upvarepsilon}(\uptheta ) } .
 \end{align}
 Experiment indicates that the limit of $j_{\upvarepsilon}(\uptheta )$ does not exist as $\upvarepsilon\rightarrow 0$ (except for $\uptheta\in\Q$, see Proposition \ref{ratsing} below), but 
 instead gives rise to a set of limit points.
 We indicate this state of affairs by writing
 \[  j^{\rm qt}(\uptheta ) := \underset{\upvarepsilon\rightarrow 0}{\text{\rm lim-pnt }}  j_{\upvarepsilon}(\uptheta )\subset\R\cup\{\infty\} \]
 where by $\text{\rm lim-pnt}_{\upvarepsilon\rightarrow 0}$ we mean the set of all limits of convergent sequences $\{  j_{\upvarepsilon_{i}}(\uptheta )\}$,
 $i=1,2,\dots$, 
 where $\upvarepsilon_{i}\rightarrow 0$.   In this way we obtain a {\it multivalued map}
 \[ j^{\rm qt}:\R\multimap \R\cup\{\infty\}, \]
 whose values may be thought of as being the spectrum of some linear operator.

 There are two privileged limit points:  
 \begin{enumerate}
 \item $ \underline{j}^{\rm qt}(\uptheta )=\liminf_{\upvarepsilon\rightarrow 0}  j^{\rm qt}_{\upvarepsilon}(\uptheta )$.
 \item $ \overline{j}^{\rm qt}(\uptheta )=\limsup_{\upvarepsilon\rightarrow 0}  j^{\rm qt}_{\upvarepsilon}(\uptheta )$.
 \end{enumerate}
 
 There is also a privileged submultimap.  Consider the sequence $\{ N_{i}\}\subset\N$ of best approximations \cite{Ca}
 to $\uptheta$
 and define $\{ \upvarepsilon_{i}^{\rm best}\}$
by
 \[  \|N_{i}\uptheta\|=\upvarepsilon_{i}^{\rm best} .\]
Then we define 
\[ j^{\rm qt}_{\sf {\small  best}}:\R\multimap \R\cup\{\infty\}\]
to be the set of limit points of  $\big\{    j^{\rm qt}_{ \upvarepsilon_{i}^{\rm best}}(\uptheta )\big\}$.
 We will see in the sequel that in the case of $\uptheta=\upvarphi=$ the golden mean, $j^{\rm qt}_{\sf {\small  best}}(\upvarphi )$ is a singleton, and experimentally satisfies (see the Appendix)
 \[  j^{\rm qt}_{\sf {\small  best}}(\upvarphi )= \underline{j}^{\rm qt}(\upvarphi ).\]

 Recall that the projective general linear group ${\rm PGL}_{2}(\Z )$ acts on $\R-\Q$ by M\"{o}bius transformations.
 
 \begin{theo} $j^{\rm qt}$ and $j^{\rm qt}_{\sf {\small  best}}$ are modular invariants: 
 \[ j^{\rm qt}(A(\uptheta ))=j^{\rm qt}(\uptheta)\quad \text{and}\quad
 j^{\rm qt}_{\sf {\small  best}}(A(\uptheta ))=j^{\rm qt}_{\sf {\small  best}}(\uptheta)\] for all $A\in {\rm PGL}_{2}(\Z )$.
\end{theo}

 \begin{proof}  Let  \[ A=\left(\begin{array}{cc}
a & b \\
c & d \\
\end{array}
\right)\] be a (representative of an) element of  ${\rm PGL}_{2}(\Z )$.  We claim that for $\upvarepsilon >0$ sufficiently small
the map
\[ n\mapsto cn^{\perp}+dn \]
defines a bijection $ B_{\upvarepsilon}(\uptheta )\leftrightarrow  B_{\upvarepsilon'}(A(\uptheta ) )$
where $\upvarepsilon' = \upvarepsilon\cdot |c\uptheta +d|^{-1}$.  Indeed
\begin{align*} A(\uptheta )(cn^{\perp}+dn) 
 &=  \frac{a\uptheta +b}{c\uptheta +d}\cdot \big( n(c\uptheta + d) -c\upvarepsilon(n) \big) \\
 & = a\uptheta n+bn -c\upvarepsilon (n)\cdot  \frac{a\uptheta +b}{c\uptheta +d} \\
 & = an^{\perp}  + bn + a\upvarepsilon (n) - c\upvarepsilon (n)\cdot  \frac{a\uptheta +b}{c\uptheta +d}  \\
 & = an^{\perp}  + bn + \frac{\upvarepsilon (n)}{c\uptheta +d}.
\end{align*}
Therefore, if $\{ \upvarepsilon_{i}\}$ produces a limit point of $j^{\rm qt}(\uptheta)$, we have
\begin{align*}   \lim_{i\rightarrow \infty}J_{\upvarepsilon_{i}}(A(\uptheta ) ) = & \lim_{i\rightarrow \infty} \frac{49}{40} \frac{\big( \sum_{n\in B_{\upvarepsilon_{i}}(\uptheta )}  (cn^{\perp} + dn)^{-6} \big)^{2}   }{\big( \sum_{n\in B_{\upvarepsilon_{i}}(\uptheta )}   (cn^{\perp} + dn)^{-4} \big)^{3}}  \\
  & \\
= & \lim_{i\rightarrow \infty}  \frac{49}{40} 
 \frac{\big( \sum_{n\in B_{\upvarepsilon_{i}}(\uptheta )}  n^{-6}(c(n^{\perp}/n)+d)^{-6} \big)^{2}   }{\big( \sum_{n\in B_{\upvarepsilon_{i}}(\uptheta )}  n^{-4} (c(n^{\perp}/n)+d)^{-4} \big)^{3}} \\
  & \\
  & = \lim_{i\rightarrow \infty} \frac{49}{40} 
 \frac{\big( \sum_{n\in B_{\upvarepsilon_{i}}(\uptheta )}  n^{-6} \big)^{2}   }{\big( \sum_{n\in B_{\upvarepsilon_{i}}(\uptheta )}  n^{-4} \big)^{3}} = 
  \lim_{i\rightarrow \infty}J_{\upvarepsilon_{i}}(\uptheta  ),\\
\end{align*}
giving the modularity of $j^{\rm qt}$.  Since
$A\in {\rm GL}_{2}(\Z )$ takes tails of best approximations to tails of best approximations (see \cite{Ca}, page 9), a similar argument
gives the modularity of $j^{\rm qt}_{\sf {\small  best}}$.
 \end{proof}
 
 We call $j^{\rm qt}(\uptheta )$ the {\bf {\small quantum modular invariant}} of $\uptheta$; by the above result, $j^{\rm qt}$
 defines a non continuous multivalued function
 \[ j^{\rm qt}: {\rm PGL}_{2}(\Z )\backslash  \overline{\R} \multimap \overline{\R} \]
 where $\overline{\R}=\R\cup \{\infty\}$.
 
 \begin{prop}\label{ratsing}  $j^{\rm qt}(\uptheta )=\infty$ for all $\uptheta\in \Q$.
 \end{prop}

 \begin{proof}     If $\uptheta = q=a/b$ written in lowest terms then for $\upvarepsilon$ sufficiently
 small, $B_{\upvarepsilon}(q)=(b)$.  For such $\upvarepsilon$,
 \[ J_{\upvarepsilon}(q) =\frac{49}{40} \frac{\big( \sum_{n\in (b), n>0}  n^{-6} \big)^{2}   }{\big( \sum_{n\in (b),n>0}   n^{-4} \big)^{3}}  = \frac{49}{40}\frac{\zeta (6)^{2}}{\zeta(4)^{3}} =1 \]
 where the last equality follows from the Euler identities $\upzeta (4) = \uppi^{4}/2\cdot 3^{2}\cdot 5$ and 
  $\upzeta (6) = \uppi^{6}/3^{3}\cdot 5\cdot 7$.
 \end{proof}

That $j^{\rm qt}\equiv\infty$ on the rationals is in keeping
with the fact that the orbit of $\Q$ is regarded as the ideal boundary of the moduli space of quantum tori
\[ {\sf Mod}^{\rm qt}=  {\rm PGL}_{2}(\Z )\backslash ( \R- \Q). \]

The next two sections are devoted to calculating $j^{\rm qt}_{\sf {\small  best}}$ of the golden mean, showing that it has, in particular, a single finite value.

\section{Modular Invariant of the Golden Mean I: An Explicit Formula}\label{golden}

Let
\[  \upvarphi := \frac{1+\sqrt{5}}{2}\]
be the golden mean.
In this section we will produce, assuming that $j^{\rm qt}_{\sf {\small  best}}(\upvarphi )$ converges, an explicit formula for $j^{\rm qt}_{\sf {\small  best}}(\upvarphi )$ obtained by evaluating at $\upvarphi $ a certain rational
expression involving weighted variants of the Rogers-Ramanujan functions.  The convergence of $j^{\rm qt}_{\sf {\small  best}}(\upvarphi )$ 
will be proved in \S \ref{golden2}; in the Appendix we will present evidence that suggests that its value is the minimum of $j^{\rm qt}(\upvarphi )$.  

We begin by recalling some facts about the golden mean and its diophantine approximations, see
for example \cite{Sch}, \cite{Vo}. 
The minimal polynomial of $\upvarphi$ is $X^{2}-X-1$
and $\upvarphi$ is a unit in $\Q (\sqrt{5})$, whose inverse is $-1$ times its conjugate:
\[   \upvarphi^{-1} = -\upvarphi' =\frac{\sqrt{5}-1}{2}.   \]
The discriminant of $\upvarphi$ is $\sqrt{5}$, and the class number of $\Q (\sqrt{5})$ is one.  
The pseudo lattice $\langle 1,\upvarphi\rangle$ has endomorphism ring equal to $O_{K}$, hence
has conductor $f=1$. 

If we denote by $[a_{0},a_{1},\dots ]$ the sequence of partial quotients of a real number $\uptheta$ then for
$\uptheta=\upvarphi$, $a_{i}=1$ for all $i$.
It follows that the sequence of best approximations $(p_{m},q_{m})$ of $\upvarphi$ is given by $(F_{m+1}, F_{m})$, where
 $\{F_{m}\}=\{1,1, 2,3,5,8,\dots \}$, $m\geq 1$, denotes the Fibonacci sequence:
\[ F_{m+1}= F_{m}+ F_{m-1}, \quad m\geq 1. \]
See for example \cite{Sch}.
This means that as $m\rightarrow\infty$,
\begin{equation}\label{errorterm}  \upvarepsilon_{m} := F_{m}\upvarphi -F_{m+1}\longrightarrow 0 
\end{equation}
 and that for all $0<n<F_{m}$, 
 \[  \|n\upvarphi \| > \| F_{m}\upvarphi  \| = | \upvarepsilon_{m}|, \]
 where as before $\|x\|$ is the distance of $x$ to the nearest integer.
 
 We recall Binet's formula \cite{NZM}:
 \[  F_{m}= \frac{\upvarphi^{m}-(\upvarphi')^{m}}{\sqrt{5}} = \frac{\upvarphi^{m}-(-1)^{m}\upvarphi^{-m}}{\sqrt{5}} 
 =\left\{\begin{array}{cc}
          \frac{\upvarphi^{m}-\upvarphi^{-m}}{\sqrt{5}}  & \text{if $m$ is even} \\
            \frac{\upvarphi^{m}+\upvarphi^{-m}}{\sqrt{5}}  & \text{if $m$ is odd.}
            \end{array}\right.
                        \]
 Using Binet's formula, we may obtain the following explicit expression for
 $\upvarepsilon_{m}$ of (\ref{errorterm}):
 \begin{align}\label{explicitformoferror} \upvarepsilon_{m}=(-1)^{m+1}\upvarphi^{-m}.
 \end{align}
 Indeed, for each integer $m$ we have
 \begin{align*}F_m\upvarphi - F_{m + 1}  & = 
 \left(\frac{\upvarphi^m + (-1)^{m + 1}\upvarphi^{-m}}{\sqrt{5}}\right)\upvarphi - \left(\frac{\upvarphi^{m + 1} + (-1)^{m}\upvarphi^{-m - 1}}{\sqrt{5}} \right)\\
 &  = 
 \frac{1}{\sqrt{5}}\left(\upvarphi^{m + 1} + (-1)^{m + 1}\upvarphi^{-m +  1} - \upvarphi^{m + 1} + (-1)^{m+1}\upvarphi^{-m -  1}\right)  \\
&  = 
 \frac{1}{\sqrt{5}}(\upvarphi + \upvarphi^{-1})(-1)^{m + 1}\upvarphi^{-m} = (-1)^{m + 1}\upvarphi^{-m}.
\end{align*}
Notice then that for $m\geq 2$, we have
\[  \|F_{m}\upvarphi\| = |\upvarepsilon_{m}| \]
and in particular,
\[ \upvarepsilon_{m}^{\rm best}=\upvarphi^{-m}.\]
 For $m$ large, $\sqrt{5}F_{m}\approx\upvarphi^{m}$, with an error term = $\pm\upvarphi^{-m}$ that decays exponentially
 as $m\rightarrow\infty$.  

Finally, we recall Zeckendorf's representation (which is actually a special case of a more general
result of Ostrowski \cite{Os}):
\begin{theo}[Zeckendorf, \cite{Ze}]\label{Zeck}  Every natural number $n\in\N$ may be written uniquely as a sum of non-consecutive
Fibonacci numbers:
\[  n = F_{I} := F_{i_{1}}+\cdots + F_{i_{k}},\quad 2\leq i_{1}, i_{1}+2\leq i_{2},\dots , i_{k-1}+2\leq i_{k},\;\; 1\leq k.    \]
\end{theo}

\begin{note}  The condition that $i_{1}\geq 2$ is to ensure uniqueness in the decomposition, otherwise the value $1$ could occur in two different ways, as $F_{1}$ or $F_{2}$.
\end{note}

We now develop an explicit formula for $j^{\rm qt}_{\sf {\small  best}}(\upvarphi)$.
Write $\upvarepsilon=|\upvarepsilon_{m}|$ and
$ B=B_{m}(\upvarphi )=\big\{ n\in\N\, |\; \| n\upvarphi \|< \upvarepsilon \big\}  $
so that
\[ J^{\rm qt}_{\upvarepsilon}(\upvarphi ) =  \frac{49}{40} \frac{\big( \sum_{n\in B}  n^{-6} \big)^{2}   }{\big( \sum_{n\in B}  n^{-4} \big)^{3}}. \]
The first step is to determine the elements of $B$ in terms of their Zeckendorf representations.
In what follows, for a multi-index $I=(i_{1},\dots, i_{k})$, define $|I|=k$.

\begin{lemm}\label{ZeckLemm}  Let $n =F_{I}=F_{i_{1}}+\cdots + F_{i_{k}}$ written in its unique Zeckendorf form.  Then
$n\in B$ if and only 
\begin{itemize}
\item[I.]  $|I|\geq 1$, $i_{1}\geq m+1$ or
\item[II.] $|I|\geq 2$, $i_{1}=m$ and $ i_{2}-m$ is odd.
\end{itemize}
\end{lemm}

\begin{note} Since the Zeckendorf form consists of sums of nonconsecutive 
Fibonacci numbers, we must have that  $ i_{2}-m\geq 3$
in II. 
\end{note}

\begin{proof}  First note that we have trivially by (\ref{explicitformoferror}) that $F_{m+i}\in B$ for $i\geq 1$. Suppose that $n=F_{I}$ is a sum of more than one non-consecutive Fibonacci numbers
and $i_{1}\geq m+1$.  Then we have
\[     \|n\upvarphi\| < \upvarphi^{-(m+1)} + \upvarphi^{-(m+3)}+\cdots = \upvarphi^{-(m+1)}(1-\upvarphi^{-2})^{-1}. \]
Since $\upvarphi =\upvarphi^{2}-1$ it follows that $(1-\upvarphi^{-2})^{-1}=\upvarphi$.   Then $ \|n\upvarphi\|<\upvarphi^{-m}$ which
implies that $n\in B$.  Thus every element of the type described in I. belongs to $B$.
On the other hand, if $i_{1}\leq m-1$, then we claim that 
\[ \upvarphi^{-m} =\upvarepsilon< \|n\upvarphi\|< 1-\upvarphi^{-1} = \upvarphi^{-2}.\] 
Indeed, if $n=F_{I}$, the associated error term sum 
\[ \upvarepsilon_{I} := \pm \upvarepsilon_{i_{1}} \pm \cdots \pm \upvarepsilon_{i_{k}}\]
is minimized in absolute value by taking $i_{1}=m-1$ and assuming
that the remaining indices $i_{2},\dots$ are such that the signs of the associated error terms $\upvarepsilon_{i_{2}},\dots $ are different from the sign
of the error term $\upvarepsilon_{m-1}$.  More precisely,
\[   |\upvarepsilon_{I}| > \upvarphi^{-(m-1)} - (\upvarphi^{-(m+2)} + \upvarphi^{-(m+4)} +\cdots ) = 
\upvarphi^{-m}(\upvarphi -\upvarphi^{-2}(1-\upvarphi^{-2})^{-1}) . \]
Since $\upvarphi^{-2}(1-\upvarphi^{-2})^{-1}=\upvarphi^{-1}$ and $\upvarphi-\upvarphi^{-1}=1$, it follows that
$ |\upvarepsilon_{I}|>\upvarphi^{-m} =\upvarepsilon$.
In addition $|\upvarepsilon_{I}|$ is maximized by taking $i_{1}=2$, $i_{2}=4,\dots$, so that
\[ |\upvarepsilon_{I}| < \upvarphi^{-2} + \upvarphi^{-4}+\cdots = \frac{1}{\upvarphi^{2}-1} = \upvarphi^{-1}.  \]
Note that the distance of the latter bound $\upvarphi^{-1}$ to the nearest integer is $1-\upvarphi^{-1}=\upvarphi^{-2}$.
It follows then from the definition of $\| \cdot \|$ and the fact that we are assuming that $m>2$ that
$ \|n\upvarphi\|>\upvarphi^{-m} =\upvarepsilon$ and
$n\not\in B$.
Now if $i_{1}=m$ and $i_{2}-m$ is even, then the error terms $\upvarepsilon_{m}$ and
$\upvarepsilon_{i_{2}}$ share the same sign, and we have 
\[   \|n\upvarphi\|> \upvarphi^{-m} + \upvarphi^{-i_{2}} - \left(\upvarphi^{-(i_{2} +3)} + \upvarphi^{-(i_{2} +5)} + \cdots \right) 
= \upvarepsilon + (\upvarphi^{-i_{2}} - \upvarphi^{-(i_{2} +3)}(1-\upvarphi^{-2})^{-1}) >\upvarepsilon 
\]
Indeed, the last inequality follows since 
\begin{align*} 
 & \upvarphi^{-i_{2}} - \upvarphi^{-(i_{2} + 3)}(1-\upvarphi^{-2})^{-1} =  \\
 & \upvarphi^{-i_{2}}(1 - \upvarphi^{-3}(1 - \upvarphi^{-2})^{-1}) = \\
 & \upvarphi^{-i_{2}}(1 - \upvarphi^{-2}(\upvarphi - \upvarphi^{-1})^{-1}) = \\
 & \upvarphi^{-i_{2}}(1-\upvarphi^{-2}) >0.
\end{align*}
On the other hand, if $i_{1}=m$ and $i_{2}=m+k$, $k$ odd, then
the sign of the corresponding error terms differ, and we have
\begin{align*}  \|n\upvarphi\| & < \upvarphi^{-m}-\upvarphi^{-m-k}+\upvarphi^{-m-k-3}+\upvarphi^{-m-k-5}+\cdots  \\
 & = \upvarphi^{-m} - \upvarphi^{-m-k}\left( 1- (\upvarphi^{-3}+\upvarphi^{-5}+\cdots )   \right) \\
 & =  \upvarphi^{-m} - \upvarphi^{-m-k}\left( 1-\upvarphi^{-3}(1-\upvarphi^{-2})^{-1}\right) \\
& =  \upvarphi^{-m} - \upvarphi^{-m-k}\left( 1-\upvarphi^{-2}\right) <\upvarepsilon 
\end{align*}
so that $n\in B$.
\end{proof}

    Let $\mathfrak{N} $ be the set of increasing, finite tuples $I=(i_{1},\dots , i_{l})$ of natural numbers with $|I|=l\geq 2$ and
which are not consecutive i.e. $ i_{1}+2\leq i_{2}, \dots , i_{l-1}+2\leq i_{l}$.  Denote
by 
\begin{align}\label{defNm}
 \mathfrak{N}(m) & = \{ I=(i_{1},\dots , i_{l})\in \mathfrak{N}|\; i_{1}\geq m\}.
 \end{align} 
Also denote by 
\begin{align}\label{defMm} 
\mathfrak{M}(m) & = \{ I\in \mathfrak{N}(m)|\; i_{1}=m \text{ and } i_{2}=m+k \text{ for }
k \text{ odd}\}.
\end{align}

Consider $B_{m}$ for $m>2$.  Then by the Lemma we have  
\[     J^{\rm qt}_{\upvarepsilon_{m}} (\upvarphi )  = \frac{49}{40}\frac{ \bigg(  \sum_{n\in B_{m}} n^{-6}  \bigg)^{2}  }{\bigg(\sum_{n\in B_{m}} n^{-4}  \bigg)^{3}} = \frac{49}{40}\frac{  \bigg( \sum_{i=1}^{\infty}  F_{m+i}^{-6}  +   \sum_{I\in\mathfrak{N}(m+1)}  F_{I}^{-6}  + 
 \sum_{ I\in \mathfrak{M}(m)} F_{I}^{-6} \bigg)^{2} }
 {  \bigg( \sum_{i=1}^{\infty}  F_{m+i}^{-4}  +   \sum_{I\in\mathfrak{N}(m+1)}  F_{I}^{-4}  + 
 \sum_{ I\in \mathfrak{M}(m)} F_{I}^{-4} \bigg)^{3}}  , \]
 an expression whose status is still only formal.  Consider also the formal expression
 \begin{align}\label{limitoftheJs}
  J_{0}^{\rm qt}(\upvarphi ) & := \frac{49}{40}\frac{ \bigg( \sum_{i= 1}^{\infty}  \upvarphi^{-6i}  +   \sum_{I\in\mathfrak{N}(1)}   \upvarphi_{I}  ^{-6}  + 
 \sum_{ I\in \mathfrak{M}(0)} \upvarphi_{I}^{-6} \bigg)^{2} }
 { \bigg( \sum_{i=1}^{\infty} \upvarphi^{-4i}   +   \sum_{I\in\mathfrak{N}(1)}   \upvarphi_{I}  ^{-4}  + 
 \sum_{ I\in \mathfrak{M}(0)} \upvarphi_{I}^{-4} \bigg)^{3}}    
 \end{align}
 where 
 \[   \upvarphi_{I} := \upvarphi^{i_{1}} +\cdots +  \upvarphi^{i_{l}} .\]

\begin{theo}\label{limitgoldenmean} If $J_{0}^{\rm qt}(\upvarphi )$ converges then so does $J^{\rm qt}_{\upvarepsilon_{m}}(\upvarphi )$ for each $m$ and
 \[    J^{\rm qt}_{\upvarepsilon_{m}}(\upvarphi ) \longrightarrow J_{0}^{\rm qt}(\upvarphi )=J_{\sf best}^{\rm qt}(\upvarphi )\]
 as $m\rightarrow\infty$.  
 \end{theo}

 \begin{proof}
 Multiply the numerator and denominator of $J^{\rm qt}_{\upvarepsilon_{m}}(\upvarphi )$ by $F^{12}_{m}$
to obtain
\begin{align}\label{FibNormOfJ} J^{\rm qt}_{\upvarepsilon_{m}}(\upvarphi ) & = \frac{49}{40}\frac{  \bigg( \sum_{i=1}^{\infty}  (F_{m}/F_{m+i})^{6}  +   \sum_{I\in\mathfrak{N}(m+1)}  (F_{m}/F_{I})^{6}  + 
 \sum_{ I\in \mathfrak{M}(m)} (F_{m}/F_{I})^{6} \bigg)^{2} }
 {  \bigg( \sum_{i=1}^{\infty}  (F_{m}/F_{m+i})^{4}  +   \sum_{I\in\mathfrak{N}(m+1)}  (F_{m}/F_{I})^{4}  + 
 \sum_{ I\in \mathfrak{M}(m)} (F_{m}/F_{I})^{4}\bigg)^{3}}  .  
 \end{align}

It will suffice to show that each term $T^{-6}_{m}=T^{-6}_{m,I}$  ($T^{-4}_{m}=T^{-4}_{m,I}$) appearing in a sum contained in the numerator (denominator) of (\ref{FibNormOfJ})
satisfies 
\[ C_{m}^{-6} \cdot T^{-6} < T^{-6}_{m} <C_{m}^{6} \cdot T^{-6}  \quad \bigg(C_{m}^{-4} \cdot T^{-4} < T^{-4}_{m} <C_{m}^{4} \cdot T^{-4}  \bigg)   \]
where $T=T_{I}$ is the correspondingly indexed term of $J^{\rm qt}_{0}(\upvarphi )$ and
\[  C_{m} = \frac{1+\upvarphi^{-2m}}{1-\upvarphi^{-2m}}.\]
This will give convergence of each $J^{\rm qt}_{\upvarepsilon_{m}}(\upvarphi )$, as well as the bound
\[   \left( \frac{1-\upvarphi^{-2m}}{1+\upvarphi^{-2m}}\right)^{24}J^{\rm qt}_{0}(\upvarphi )<  J^{\rm qt}_{\upvarepsilon_{m}}(\upvarphi ) < \left( \frac{1+\upvarphi^{-2m}}{1-\upvarphi^{-2m}}\right)^{24}J^{\rm qt}_{0}(\upvarphi ) ,  \]
which implies that $J^{\rm qt}_{\upvarepsilon_{m}}(\upvarphi ) \rightarrow J^{\rm qt}_{0}(\upvarphi )$.

We will now make use of Binet's formula, $\sqrt{5}F_{m}=(\upvarphi^{m}\pm\upvarphi^{-m})$.  
Note that
the $\sqrt{5}$ factors drop out and so we may simply replace every Fibonacci term $F_{m}$ appearing by
$\upvarphi^{m}\pm\upvarphi^{-m}$.

We consider first the numerator of (\ref{FibNormOfJ}), treating each of the three sums there separately.  The
first sum may be written \[   \sum_{i=1}^{\infty}  (F_{m}/F_{m+i})^{6} = \sum_{i=1}^{\infty}  \left( \frac{\upvarphi^{m} \pm\upvarphi^{-m}}{ \upvarphi^{m+i} \pm(-1)^{i} \upvarphi^{-(m+i)}}  \right)^{6} 
 =\sum_{i=1}^{\infty} \upvarphi^{-6i}  \left( \frac{1 \pm\upvarphi^{-2m}}{1 \pm(-1)^{i}  \upvarphi^{-2m-2i}}  \right)^{6} .        \]
Note that
 \[  \left( \frac{1-\upvarphi^{-2m}}{1+\upvarphi^{-2m}}\right)^{6} <   \left( \frac{1 \pm\upvarphi^{-2m}}{1 \pm(-1)^{i}  \upvarphi^{-2m-2i}}  \right)^{6} <
\left( \frac{1+\upvarphi^{-2m}}{1-\upvarphi^{-2m}}\right)^{6}  .  \]
The next sum is 
\begin{align}\label{BinetNormOfJ}
     \sum_{I\in\mathfrak{N}(m+1)}   \left( F_{m}/ F_{I}\right)^{6} & = \sum_{I\in\mathfrak{N}(m+1)}   \left(\frac{ \upvarphi^{m}\pm\upvarphi^{-m} }{ (\upvarphi^{m+i_{1}}\pm \upvarphi^{-m-i_{1}} )+\cdots + (\upvarphi^{m+i_{k}}\pm\upvarphi^{-m-i_{k}})  } \right)^{6}, 
     \end{align}
where we are writing our generic $I\in\mathfrak{N}(m+1)$ in the form $I=(i_{1}+m,\dots ,i_{k}+m)$ with $1\leq i_{1}<i_{2}<\cdots <i_{k}$. 
Letting $I_{0}=(i_{1},\dots , i_{k})$ then each term of the sum in (\ref{BinetNormOfJ}) may be re-written
\begin{align}\label{rewritingtheterms}  \left( \frac{1\pm\upvarphi^{-2m}}{ \upvarphi_{I_{0}} + (\pm\upvarphi_{-I_{0}-2m}) }\right)^{6}
=   \upvarphi_{I_{0}}^{-6} \cdot\left( \frac{1\pm\upvarphi^{-2m}}{1 + (\pm \upvarphi_{-I_{0}-2m})/  \upvarphi_{I_{0}}  }\right)^{6}
 \end{align}
 where 
 \[\pm\upvarphi_{-I_{0}-2m} :=\pm \upvarphi^{-i_{1}-2m}\pm \cdots\pm  \upvarphi^{-i_{k}-2m}, \]
 the signs determined as in Binet's formula by the parities of the powers.
It is easy to see that
 \begin{align}\label{basicinequalities}
   \left( \frac{1-\upvarphi^{-2m}}{1+\upvarphi^{-2m}}\right)^{6}
   < \left( \frac{1\pm\upvarphi^{-2m}}{1 + (\pm \upvarphi_{-I_{0}-2m})/  \upvarphi_{I_{0}}  }\right)^{6}< \left( \frac{1+\upvarphi^{-2m}}{1-\upvarphi^{-2m}}\right)^{6}:
   \end{align}
indeed, both inequalities in (\ref{basicinequalities}) follow since
 \[  \upvarphi^{-2m} > (\pm \upvarphi_{-I_{0}-2m})/  \upvarphi_{I_{0}} >-\upvarphi^{-2m},\]
 true as
 \begin{align}\label{trueas}  (\pm \upvarphi_{-I_{0}-2m})/\upvarphi_{I_{0}} = \upvarphi^{-2m}
\left( \frac{\pm\upvarphi^{-i_{1}}\pm \cdots\pm  \upvarphi^{-i_{k}} }{\upvarphi^{i_{1}}+ \cdots+  \upvarphi^{i_{k}}}\right).
  \end{align}
 What remains is the sum over $\mathfrak{M}(m)$: the analysis here is essentially the same as that made
 for the sum over $\mathfrak{N}(m+1)$, only we take into account that $I=(m,m+j,m+i_{3},\dots ,m+i_{k})$
 where $j$ is odd.  Writing $I_{0}=(0,j,i_{3},\dots ,i_{k})$, then we have the equation (\ref{rewritingtheterms}) with
 \[  \pm\upvarphi_{-I_{0}-2m} =\pm\upvarphi^{-2m}\mp\upvarphi^{-j-2m} \pm\cdots\pm  \upvarphi^{-i_{k}-2m}, \]
where the $\mp$ sign of $\upvarphi^{-j-2m}$ indicates that this sign is opposite to that of $\upvarphi^{-2m}$, as  
 $j$ is odd.  The analogue of (\ref{trueas}) is then
 \begin{align*} (\pm \upvarphi_{-I_{0}-2m})/\upvarphi_{I_{0}} = \upvarphi^{-2m}
\left( \frac{\pm1\mp\upvarphi^{-j}\pm\cdots\pm  \upvarphi^{-i_{k}} }{1 +\upvarphi^{j}+ \cdots+  \upvarphi^{i_{k}}}\right),
  \end{align*}
  which yields the analogue of (\ref{basicinequalities}) in this case.  This completes our bounding of
  the numerator.  Analogous bounds, with the exponent $6$ replaced by $4$, may be found for the corresponding sums in the denominator of $J^{\rm qt}_{\upvarepsilon_{m}}$.
 The result now follows.
 \end{proof}
 
Let $P(n)$ be the set of partitions of $n$ into
  into distinct parts whose differences are at least $2$, and let $c(n)=|P(n)|$. 
   The generating function
  \[  F(x)=\sum_{n=1}^{\infty} c(n)x^{n}=\sum \frac{x^{n^{2}}}{(1-x)\cdots (1-x^{n})}\] 
 is of substantial combinatorial interest: $1+F(x)$ is the left-hand side of the first Rogers-Ramanujan identity \cite{HaWr}.
 
For each partition $I\in P(n)$, let $f_{I}(x)=x^{i_{1}}+\cdots +x^{i_{k}}$
be the associated weighting polynomial.
Define
 \[ C_{x,M}(n)=x^{Mn} \sum _{I\in P(n) }
  f_{I}(x)^{-M}.\]
  Consider the generating function
  \[  G_{M}(x)=\sum C_{x,M}(n) x^{n}.\]
  Clearly we have 
  \[   G_{M}(\upvarphi) =  \sum_{i= 1}^{\infty}  \upvarphi^{-Mi}  +   \sum_{I\in\mathfrak{N}(1)}   \upvarphi_{I}  ^{-M} . \]
  Similarly, let $Q(n)\subset P(n)$ 
be the set of those partitions $I=i_{1}<i_{2}<\cdots <i_{k}$ in $P(n)$ for which $i_{1}$ is odd and $\geq 3$.
Let 
\[ D_{x,M}(n) :=x^{Mn}  \sum _{I\in Q(n)}(1+f_{I}(x ))^{-M}\]
and define \[  H_{M}(x):=\sum D_{x,M}(n) x^{n}.\]
Then 
 \[   H_{M}(\upvarphi) =  \sum_{ I\in \mathfrak{M}(0)} \upvarphi_{I}^{-M}.\]
 The following is then immediate:

\begin{coro}\label{explicitformula}  Let $J_{0}^{\rm qt}(\upvarphi )$ be as above.  Then
\begin{align}\label{rogramexpress} J_{0}^{\rm qt}(\upvarphi )=\frac{49}{40}\frac{ \big(G_{6}(\upvarphi)  + H_{6}(\upvarphi)\big)^{2}}{ \big(G_{4}(\upvarphi)  + H_{4}(\upvarphi) \big)^{3}} .
\end{align} 
\end{coro}

 \begin{note} 
  If one replaces in the formula for $C_{x, M}(n)$ the weighting polynomial $ f_{I}(x )^{-M}$ by the equiweight $x^{-Mn}$ one recovers $c(n)$.  Thus the functions $G_{M}(x), H_{M}(x)$ may be viewed as weighted variants of the variable part of the Rogers-Ramanujan function.
  \end{note}

In the Appendix, we show that by replacing $G_{M}(\upvarphi)$ by $G_{M}'(\upvarphi)= G_{M}(\upvarphi)+1$ in (\ref{rogramexpress}) we obtain a value
very close to the experimental supremum of $J^{\rm qt}(\upvarphi )$.

 \section{Modular Invariant of the Golden Mean II: Convergence}\label{golden2}
 
In this section we will show that $j^{\rm qt}_{\sf best}(\upvarphi )<\infty$ and in fact converges.  As before we write
 $   j^{\rm qt}_{\sf best}(\upvarphi ) := 12^{3}/( 1-J_{\sf best}^{\rm qt}(\upvarphi )  ) $.  
 
\begin{theo}\label{goldenbound} $j_{\sf best}^{\rm qt}(\upvarphi )$ converges with the bounds 
 \[   9150< j^{\rm qt}_{\sf best}(\upvarphi ) <9840.  \]
 \end{theo}

 \begin{proof}  To prove the convergence of
 $j_{\sf best}^{\rm qt}(\upvarphi )$, it is enough to prove convergence of the explicit formula $j^{\rm qt}_{0}(\upvarphi ):=12^{3}/( 1-J_{\sf best}^{\rm qt}(\upvarphi )  ) $
 obtained from (\ref{limitoftheJs}).
 Observe first that
 \[ \sum_{i= 1}^{\infty}  \upvarphi^{-6i} = (\upvarphi^{6}-1)^{-1},\quad\quad  \sum_{i=1}^{\infty}  \upvarphi^{-4i} = (\upvarphi^{4}-1)^{-1}  \]
 so we may write
  \[    J_{0}^{\rm qt}(\upvarphi ) = \frac{49}{40}\frac{ \left( (\upvarphi^{6}-1)^{-1}  +   \sum_{I\in\mathfrak{N}(1)}   \upvarphi_{I}  ^{-6}  + 
 \sum_{ I\in \mathfrak{M}(0)} \upvarphi_{I}^{-6} \right)^{2} }
 { \left( (\upvarphi^{4}-1)^{-1} + \sum_{I\in\mathfrak{N}(1)}   \upvarphi_{I}  ^{-4}  + 
 \sum_{ I\in \mathfrak{M}(0)} \upvarphi_{I}^{-4} \right)^{3}}  .  \]
 We now find an explicit approximation and an upper bound for the sum  $\sum_{I\in\mathfrak{N}(1)}   \upvarphi_{I}^{-M} $
 where $M$ is a positive integer.  
 In fact, we will show that
 \begin{align}\label{goldennumerator1}
  \sum_{I\in\mathfrak{N}(1)}   \upvarphi_{I}^{-M}  = \frac{1}{(\upvarphi^{M}-1)(\upvarphi^{2}+1)^{M}} +C(M)  
  \end{align}
 where
 \begin{align}\label{goldennumeratorbound}
 C(M) < \widetilde{C}(M):=\frac{1}{\upvarphi^{2M}(\upvarphi^{M}-1)^{2}}   +\frac{1}{\upvarphi^{M}(\upvarphi^{M}-1)^{2}(\upvarphi^{2M}-\upvarphi^{M}-1)}. 
 \end{align}
 
 Consider first the sum of those $I$
 with $|I|=2$:
 \begin{align}
  \mathop{ \sum_{i_{1}\geq 1} }_ {i_{2}\geq i_{1}+2} \frac{1}{( \upvarphi^{i_{1}} + \upvarphi^{i_{2}}   )^{M}} & =
 \sum_{i=1}^{\infty} \upvarphi^{-Mi}\sum_{k=2}^{\infty} (1+ \upvarphi^{k})^{-M} \nonumber \\
  & = \frac{1}{\upvarphi^{M}-1}\left\{ \frac{1}{(1+\upvarphi^{2})^{M}} +\sum_{k=3}^{\infty} (1+ \upvarphi^{k})^{-M}  \right\} \label{goldenfirstorder}\\
& <  \frac{1}{\upvarphi^{M}-1}\left\{ \frac{1}{(1+\upvarphi^{2})^{M}}+\sum_{k=3}^{\infty} \upvarphi^{-Mk}  \right\}\nonumber \\
& = \frac{1}{(\upvarphi^{M}-1)(\upvarphi^{2}+1)^{M}} + \frac{1}{\upvarphi^{2M}(\upvarphi^{M}-1)^{2}} .\label{goldenfirstbound}
  \end{align}
  The equality (\ref{goldenfirstorder}) produces the explicit term 
  $1/((\upvarphi^{M}-1)(\upvarphi^{2}+1)^{M})$ appearing in (\ref{goldennumerator1}); the second term in (\ref{goldenfirstbound})
  is the first bounding term 
  in (\ref{goldennumeratorbound}).
  
  For $|I|=3$ we have
 \begin{align*} 
  \mathop{ \sum_{i_{1}\geq 1} }_ {i_{2}\geq i_{1}+2, i_{3}\geq i_{2}+2 } \frac{1}{( \upvarphi^{i_{1}} + \upvarphi^{i_{2}} + \upvarphi^{i_{3}} )^{M}}  &  
 =  \mathop{ \sum_{i_{1}\geq 1} }_ {i_{2}\geq i_{1}+2, i_{3}\geq i_{2}+2 } \upvarphi^{-Mi_{1}}\frac{1}{( 1 + \upvarphi^{i_{2}-i_{1}} + \upvarphi^{i_{3}-i_{1}} )^{M}} \\
 & <   \mathop{ \sum_{i_{1}\geq 1} }_ {i_{2}\geq i_{1}+2, i_{3}\geq i_{2}+2 } \upvarphi^{-Mi_{1}}\frac{1}{( \upvarphi^{i_{2}-i_{1}} + \upvarphi^{i_{3}-i_{1}} )^{M}} \\
 & = \mathop{ \sum_{i_{1}\geq 1} }_ {i_{2}\geq i_{1}+2, i_{3}\geq i_{2}+2 } \upvarphi^{-Mi_{1}}\upvarphi^{-M(i_{2}-i_{1})}\frac{1}{(1+ \upvarphi^{i_{3}-i_{2}} )^{M}} \\
& < \sum_{i\geq 1} \upvarphi^{-Mi}
 \sum_{j\geq 2}\upvarphi^{-Mj}\sum_{k\geq 2} \upvarphi^{-Mk} \\
 & = \frac{(\upvarphi^{-M})^{2}}{ (\upvarphi^{M}-1)^{3} }  .
 \end{align*}
 Inductively, for the terms with $|I|=l\geq 3$ we have the bound
 \[  \frac{(\upvarphi^{-M})^{l-1}}{ (\upvarphi^{M}-1)^{l} }.  \]
 Summing these bounds from $l=3$ to $\infty$ gives the second term in (\ref{goldennumeratorbound}):
 \[\sum_{l=3}^{\infty} \frac{(\upvarphi^{-M})^{l-1}}{ (\upvarphi^{M}-1)^{l} }
 = \upvarphi^{M}\sum_{l=3}^{\infty} \frac{1}{ (\upvarphi^{M}(\upvarphi^{M}-1))^{l} }
 =\frac{1}{\upvarphi^{M}(\upvarphi^{M}-1)^{2}(\upvarphi^{2M}-\upvarphi^{M}-1)}
 \]
 
 We now bound the second type of sum appearing in $J_{0}^{\rm qt}(\upvarphi)$, 
$ \sum_{ I\in \mathfrak{M}(0)} \upvarphi_{I}^{-M}  $.  We will show here that
\begin{align}
\sum_{ I\in \mathfrak{M}(0)} \upvarphi_{I}^{-M} = \frac{1}{(1+\upvarphi^{3})^{M}} + D(M)
\end{align}
where
\begin{align}
D(M) <\widetilde{D}(M):=\frac{1}{\upvarphi^{3M}(\upvarphi^{2M}-1)}  +   \frac{1}{\upvarphi^{M}(\upvarphi^{2M}-1)(\upvarphi^{2M}-\upvarphi^{M}-1)}
\end{align}

When $|I|=2$ we have, since $i_{1}=0$, that $i_{2}=2j+1$ is odd, where $j\geq 1$ (recall the definition
of $\mathfrak{M}(m)$ found in (\ref{defMm})).  For such $I$ we have the contribution
\begin{align}
  \mathop{ \sum_{i=2j+1} }_{j\geq 1}\frac{1}{(1+ \upvarphi^{i})^{M}} 
& = \frac{1}{(1+\upvarphi^{3})^{M}} + \sum_{j= 2}^{\infty} (1+\upvarphi^{(2j+1)} )^{-M} \\
& <  \frac{1}{(1+\upvarphi^{3})^{M}} +\sum_{j= 2}^{\infty} \upvarphi^{-M(2j+1)} \nonumber \\
& = \frac{1}{(1+\upvarphi^{3})^{M}} + \upvarphi^{-M}\sum_{j= 2}^{\infty}\upvarphi^{-2Mj} \nonumber \\
& = \frac{1}{(1+\upvarphi^{3})^{M}} + \upvarphi^{-5M}\frac{1}{1-\upvarphi^{-2M}} \nonumber \\
& = \frac{1}{(1+\upvarphi^{3})^{M}}+ \frac{1}{\upvarphi^{3M}(\upvarphi^{2M}-1)} .
 \end{align}

For $|I|=3$ we have
 \begin{align*} 
   \sum_ {j\geq 1, k\geq (2j+1)+2 } \frac{1}{( 1 + \upvarphi^{2j+1} + \upvarphi^{k} )^{M}}  &  
 < \sum_ {j\geq 1, k\geq (2j+1)+2 } \upvarphi^{-M(2j+1)}\frac{1}{( 1+ \upvarphi^{k-(2j+1)} )^{M}} \\
 &< \sum_{j=1}^{\infty}\upvarphi^{-M(2j+1)}\sum_{k=2}^{\infty}\upvarphi^{-Mk} \\
  & = \frac{1}{\upvarphi^{M}(\upvarphi^{2M}-1)}\cdot \frac{1}{\upvarphi^{M}(\upvarphi^{M}-1)} \\
  & = \frac{1}{\upvarphi^{M}+1}\cdot \left( \frac{\upvarphi^{-M}}{\upvarphi^{M}-1}\right)^{2} \\ 
 \end{align*}
For the sum over $I$ with $|I|=l$, we obtain inductively 
the bound
\[       \frac{1}{\upvarphi^{M}+1}\left( \frac{\upvarphi^{-M}}{\upvarphi^{M}-1}  \right)^{l-1}    \]
and summing these from $l=3$ to $\infty$ gives
\[   \frac{1}{\upvarphi^{M}(\upvarphi^{2M}-1)(\upvarphi^{2M}-\upvarphi^{M}-1)}. \]
It follows then that
\begin{align*}
  J_{0}^{\rm qt}(\upvarphi ) & <\frac{49}{40} \frac{\left( (\upvarphi^{6}-1)^{-1}  + \big( (\upvarphi^{6}-1)(\upvarphi^{2}+1)^{6}\big)^{-1}
+(1+\upvarphi^{3})^{-6} + \widetilde{C}(6) + \widetilde{D}(6) \right)^{2}}{ \left( (\upvarphi^{4}-1)^{-1}+\big( (\upvarphi^{4}-1)(\upvarphi^{2}+1)^{4}\big)^{-1}
+(1+\upvarphi^{3})^{-4}\right)^{3}} \\
& \\
& \approx 0.824376700276.
 \end{align*}
A lower bound may be given by 
\begin{align*}
  0.81115979990388 & \approx\frac{49}{40}\frac{\left( (\upvarphi^{6}-1)^{-1}  + \big( (\upvarphi^{6}-1)(\upvarphi^{2}+1)^{6}\big)^{-1}
+(1+\upvarphi^{3})^{-6} \right)^{2}}{ \left( (\upvarphi^{4}-1)^{-1}+\big( (\upvarphi^{4}-1)(\upvarphi^{2}+1)^{4}\big)^{-1}
+(1+\upvarphi^{3})^{-4}+ \widetilde{C}(4) + \widetilde{D}(4) \right)^{3}} \\
& \\
& <J_{0}^{\rm qt}(\upvarphi ) 
\end{align*}
which give the bounds presented in the statement of the theorem.  Since the numerator and denominator of $J_{0}^{\rm qt}(\upvarphi )$
are hypergeometric functions with positive coefficients evaluated at a positive real number, 
 it follows from the above bounds that they converge, and in particular, that
$J_{0}^{\rm qt}(\upvarphi )$ converges.
\end{proof}

 \begin{note}  Using the PARI/GP value of the explicit formula of Corollary \ref{explicitformula}, we get
$j^{\rm qt}_{\sf best}(\upvarphi )\approx 9538.249655644$, which agrees closely with the experimental value obtained
for $\underline{j}^{\rm qt}(\upvarphi )$.  See the Appendix.
 \end{note}

\section{Quantum Tori and Kronecker Foliations}\label{qtoriandkfol}

In this section we begin the process of finding a continuous and single valued version of the set-valued quantum invariant $j^{\rm qt}$,
as well as putting the latter in its proper geometrical context. 

Consider the Kronecker foliation $\mathcal{F}(\uptheta)$ of slope $\uptheta$ 
in the elliptic curve $\T(i) = \C/\langle 1, i\rangle$,
i.e. the image in $\T(i)$ 
of the foliation of the complex plane $\C$ by lines of slope $\uptheta$.  The leaf
space of $\mathcal{F}(\uptheta)$ may be identified with
the quotient group $\T(i) /L(\uptheta)$,
where $L(\uptheta)$ is the leaf through the origin, a 1-parameter
subgroup of $\T(i)$.
When $\uptheta\in\R -\Q$,  $L(\uptheta)$ is dense in $\T(i)$
so that the leaf space is non Hausdorff.

On the other hand, let $\Uplambda(\uptheta)=\langle 1, \uptheta \rangle\subset \R$ be the pseudo
lattice generated by $1$ and $\uptheta$.  As discussed in the Introduction, 
the quantum torus associated to $\uptheta\in\R$ may be defined
as the following quotient:
\[  \T(\uptheta) = \R/ \Uplambda(\uptheta) .\]
When $\uptheta$ is irrational, this is a non Hausdorff topological group.  
It will be convenient for us to allow $\uptheta$ to be rational as well, in which case one obtains 
the circle.

\begin{prop}\label{leafspace}  The leaf space of $\mathcal{F}(\uptheta)$ is canonically isomorphic to $\T(\uptheta)$.
\end{prop}

\begin{proof} Writing $\SI^{1}=\R/\Z$, consider the suspension 
$   (\R\times  \SI^{1} )/\Z$ ,
where the action of $\Z$ is diagonal: $n\cdot (r, s+\Z) = (r+n, (x- \uptheta n)+\Z) $.  The suspension
defines a linear foliation of $\T(i)$: the image of the product foliation $\R\times \SI^{1}$, whose leaves are of the form
$\R\times \{ s+\Z\}$.  There is an isomorphism of the Kronecker foliation $\mathcal{F}(\uptheta)$ with this foliation,  
induced by
 $\C\rightarrow (\R\times  \SI^{1} )$, $r+is\mapsto (r, s-r\uptheta+\Z)$.
Through this identification, one sees that the leaf space of $\mathcal{F}(\uptheta)$ is canonically isomorphic to the quotient group $\SI^{1}/\langle \uptheta +\Z\rangle$.  But the latter is canonically isomorphic to $\T(\uptheta)$.
\end{proof}

The Kronecker foliation has an obvious generalization in which one replaces $\T(i)$ by any elliptic curve
$ \T(\upmu ) = \C /\Uplambda(\upmu ) $
where $\Uplambda(\upmu ) =\langle 1, \upmu\rangle$ and where $\upmu\in\HP$ = the hyperbolic plane.  
Given $\uptheta\in\R\cup\{\infty\}\approx\SI^{1}$, let $\widetilde{\mathcal{F}}(\upmu ,\uptheta)$
be the foliation of $\C$ defined by the translates of the line of \guillemotleft $\upmu$-slope $\uptheta$\guillemotright ,
\[ \widetilde{L}(\upmu ,\uptheta)= \left\{ \begin{array}{ll}
                                     \R\cdot (1 + \uptheta\upmu ) & \text{if $\uptheta\not=\infty$} \\
                                       \R\cdot\upmu & \text{if $\uptheta=\infty$}
                                      \end{array}\right. \]
                                      The image 
$\mathcal{F}(\upmu ,\uptheta)$ of $\widetilde{\mathcal{F}}(\upmu ,\uptheta)$ in  $\T(\upmu )$  
is called a {\bf {\small generalized Kronecker foliation}} of slope $\uptheta$ and modulus $\upmu$.
Alternatively, $\kf$ is completely determined by the pair \[ \big(\, \T(\upmu ), L(\upmu ,\uptheta)\big)\] consisting
of the elliptic curve and the distinguished 1-parameter subgroup
\[ L(\upmu ,\uptheta) = \text{image of the line $\widetilde{L}(\upmu ,\uptheta)$} = \text{leaf
through $0$}  .\]
This may be regarded as a continuous generalization of the notion of an elliptic curve
equipped with a distinguished finite subgroup of order $N$.

As in \cite{Ma}, it will be convenient to allow the parameter $\upmu$ to take values in 
$\overline{\HP}$ as well. 
If we denote by $\pm\HP = \HP\cup\overline{\HP}$ then ${\rm PGL}_{2}(\Z )$ acts on
$\pm\HP$ by isometries and we recover by quotient the classical moduli space of elliptic curves
\[ {\rm Mod}^{\rm cl} := {\rm PSL}_{2}(\Z )\backslash \HP\approx {\rm PGL}_{2}(\Z )\backslash \pm\HP  .  \]
The Kronecker foliation 
$\mathcal{F}(\upmu ,\uptheta)$ for $\upmu\in\overline{\HP}$  is defined
exactly as in the case of $\upmu\in\HP$.
Note that for all 
$(\upmu ,\uptheta)\in \pm\HP\times\SI^{1}$ we have the equality
\[ \mathcal{F}(-\upmu , -\uptheta)=  \mathcal{F}(\upmu ,\uptheta).\]
This equality remains true for $\uptheta =\infty$ (which, like $0$, has no sign).

Let $(\upmu ,\uptheta), (\upmu' ,\uptheta ')\in \pm\HP\times\SI^{1}$.
The Kronecker foliations $\mathcal{F}(\upmu ,\uptheta)$ and
$\mathcal{F}(\upmu' ,\uptheta')$ are said to be equivalent if there exists
a homothety $z\mapsto \uplambda z$ inducing an isomorphism of underlying elliptic curves
that transports $\mathcal{F}(\upmu ,\uptheta)$ to
$\mathcal{F}(\upmu' ,\uptheta')$:  or equivalently, inducing an
isomorphism of pairs
\[ f: \bigg(\T(\upmu ), L(\upmu ,\uptheta)\bigg)\longrightarrow \bigg(\T(\upmu'), L(\upmu' ,\uptheta')\bigg) .   \]
Note that this notion of equivalence is formally in agreement with that used for pairs
of tori and finite subgroups of a fixed order $N$.

In what follows, for any $A\in {\rm PGL}_{2}(\Z )$, denote by $A^{-T}$ the contragredient class i.e.\ the transformation
defined by the inverse of the transpose of a matrix in the projective class of $A$: note that 
$(AB)^{-T}=A^{-T}B^{-T}$.

\begin{prop}  $\mathcal{F}(\upmu ,\uptheta)$ is isomorphic to
$\mathcal{F}(\upmu' ,\uptheta')$
$\Leftrightarrow$ there exists $A\in {\rm PGL}_{2}(\Z )$ such that 
\[ \upmu' = A(\upmu )\quad\text{and}\quad \uptheta' = A^{-T} (\uptheta ) . \]
\end{prop}

\begin{proof}  Assume first that $\upmu, \upmu'\in \HP$ and $\mathcal{F}(\upmu ,\uptheta)$ and
$\mathcal{F}(\upmu' ,\uptheta')$ are isomorphic via the homothety defined by $\uplambda\in\C$ with $\uplambda\cdot \Uplambda(\upmu ) =  \Uplambda (\upmu ')$.  Then we have
$\uplambda\upmu = a\upmu' + b$, $ \uplambda = c\upmu' + d $
where 
\[ B = \left( \begin{array}{ll} 
                     a & b \\
                     c & d
              \end{array}
  \right)\in {\rm SL}(2,\Z ) . \]
Thus $\upmu = B(\upmu ' )$ or $\upmu' = B^{-1}(\upmu  )$.
On the other hand,
\begin{eqnarray*} 
\uplambda \cdot \widetilde{L}(\upmu  , \uptheta)  & = &  \R\cdot \big( (c\upmu' +d) + \uptheta(a\upmu' +b)\big) \\
 & \\
& = &  \R\cdot \big( (b\uptheta +d) + (a\uptheta +c)\upmu' \big)\\
 & \\
 & = &  \R\cdot \big(1 + B^{T}(\uptheta)\upmu'\big) \\
& \\
 & = & \widetilde{L}\big(B^{-1}(\upmu  )  ,\; B^{T}(\uptheta)\big).
\end{eqnarray*}
This shows that multiplication by $\uplambda$
induces an equivalence of foliations
\[ \mathcal{F}(\upmu ,\uptheta)\longrightarrow \mathcal{F}\big(B^{-1}(\upmu  )  , B^{T}(\uptheta)\big) . \]
Writing $A=B^{-1}$ we obtain the form of equivalence stated in the Proposition.

In case $(\upmu' ,\uptheta') =(-\upmu , -\uptheta)$ then we take 
\[ A = \left( \begin{array}{rr} 
                     -1 & 0 \\
                     0 & 1
              \end{array}
  \right)\in {\rm GL}_{2} (\Z )  \]
and noting that $A=A^{T}=A^{-1}$, we have $(A(\upmu ) , A^{-T} (\uptheta )) =(-\upmu , -\uptheta)$.
The argument above is symmetric, so that if there exists $A\in {\rm PGL}_{2}(\Z )$ such that 
$\upmu' = A(\upmu )$ and $\uptheta' = A^{-T} (\uptheta )$, then the corresponding Kronecker
foliations are equivalent.
\end{proof}

The moduli space of isomorphism classes of quantum tori
is defined  \cite{Mar}
\[  {\sf Mod}^{\rm qt}=  {\rm PGL}_{2}(\Z )\backslash ( \R- \Q)\]  
which may be viewed as a kind of boundary of the classical moduli space 
${\sf Mod}^{\rm cl}$.  Note that   ${\sf Mod}^{\rm qt} $
is itself a non Hausdorff space since ${\rm PGL}_{2}(\Z )$ acts densely on $ \R-\Q$.
The orbit of $[\infty ]$ of $\Q$ is viewed as the ideal boundary of ${\sf Mod}^{\rm qt}$ and we write as well
\[ \overline{\sf Mod}^{\rm qt} :=  {\sf Mod}^{\rm qt}\cup [\infty ].  \]

By Proposition~\ref{leafspace}, the moduli space of generalized
Kronecker foliations is the \guillemotleft signed\guillemotright\ Anosov foliation
\[     {\sf Mod}^{\rm kf} =  {\rm PGL}_{2}(\Z )\backslash( \pm\HP\times \SI^{1})   \]
where $A\in {\rm PGL}_{2}(\Z )$ acts by 
\begin{align}\label{Aactiononpairs}
 A\cdot (\upmu ,\uptheta  ) = (A(\upmu ),  A^{-T} (\uptheta ) ).
 \end{align}    
We regard the images of $\pm\HP\times \{\uptheta \}$ in ${\sf Mod}^{\rm kf}$ as the leaves.
This foliation fibers over ${\sf Mod}^{\rm cl}={\rm PGL}_{2}(\Z )\backslash\pm\HP$,
in which the fiber ${\sf Mod}^{\rm kf}_{[i]}$ over the class $[i]\in {\sf Mod}^{\rm cl}$ parametrizes
the \guillemotleft classical\guillemotright\ Kronecker foliations.
We have the following moduli space analogue
of Proposition~\ref{leafspace}:

\begin{prop}
The leaf space of ${\sf Mod}^{\rm kf} $
is in canonical bijection with $\overline{\sf Mod}^{\rm qt}$. 
\end{prop}

\begin{proof} Since the fiber ${\sf Mod}^{\rm kf}_{[i]}$ is a complete transversal
of the foliation ${\sf Mod}^{\rm kf}$, the leaf space of  ${\sf Mod}^{\rm kf} $ may be identified 
with the set of leaf classes of elements of ${\sf Mod}^{\rm kf}_{[i]}$.  The latter
is the image of $\{ i \} \times \SI^{1}\subset \pm\HP\times \SI^{1}$ under the suspension
quotient.  In particular, two points of ${\sf Mod}^{\rm kf}_{[i]}$ lie on the same leaf
if and only if their preimages $(i, \uptheta )$, $(i, \uptheta' )\in \{ i \} \times \SI^{1}$ satisfy $\uptheta' = A(\uptheta )$ for some
$A\in {\rm PGL}_{2}(\Z )$ acting projective linearly on $\SI^{1}\approx \R\cup\{\infty\}$.  
Thus the leaf space of  ${\sf Mod}^{\rm kf} $ may be put in canonical
bijection with ${\rm PGL}_{2}(\Z )\backslash \SI^{1}\approx\overline{\sf Mod}^{\rm qt}$.
 
 \end{proof}

As mentioned in the Introduction, ${\sf Mod}^{\rm kf}$
provides a natural generalization of the moduli space $\Upgamma_{0}(N)\backslash\HP$
that classifies isomorphism classes of ordered pairs $(E, C)$,
where $E$ is an elliptic curve defined over $\C$,
and $C$ is a cyclic subgroup of $E$ of order $N$.  

We could extend $j^{\rm qt}$ -- which is defined on the transversal
${\sf Mod}^{\rm kf}_{[i]}$ -- to all of ${\sf Mod}^{\rm kf} $ using a similar definition to that found in \S 1, but the discontinuity
and multivaluedness would persist.  Instead we will use nonstandard models to construct an analogous space which fibers over
${\sf Mod}^{\rm kf} $ on which
 $j^{\rm qt}$ lifts to a continuous single valued function.

\section{Nonstandard Structures}\label{nonstdstruct}

In what follows, $\mathcal{I}$ is a discrete, infinite 
set.
\vspace{3mm}

\noindent  {\em Ultrafilters and Stone Spaces}

\vspace{3mm}

\noindent Recall that a filter on $\mathcal{I}$ is a subset $\mathfrak{f}\subset {\sf 2}^{\mathcal{I}}$ not containing the empty set, which is closed
with respect to finite intersections and upward inclusions ($X\in \mathfrak{f}$ and $Y\supset X$ $\Rightarrow$ $Y\in \mathfrak{f}$).
Dually the set of complements 
$I_{\mathfrak{f}} := \{ X\, | \; \mathcal{I}-X\in \mathfrak{f} \}$
is a proper ideal in the Boolean algebra ${\sf 2}^{\mathcal{I}}$.  
A maximal filter $ \mathfrak{u}$ is called an ultrafilter, whose
set of complements $I_{\mathfrak{u}}$ is a maximal ideal of ${\sf 2}^{\mathcal{I}}$.
A filter $\mathfrak{f}$ is called nonprincipal if there exists no $X\in \mathfrak{f}$ with $X\subset Y$ for all
$Y\in\mathfrak{f}$, or dually,
if $I_{\mathfrak{f}}$ is a nonprincipal ideal. See \cite{Je}.

If one has a family $\mathcal{A}\subset  {\sf 2}^{\mathcal{I}}$ of subsets not containing the empty set
and satisfying the finite intersection property, there is a unique minimal filter containing
$\mathcal{A}$, the filter $\langle \mathcal{A}\rangle$ generated by $\mathcal{A}$.
For example, if $\mathcal{I}$ is directed, $\upgamma\in \mathcal{I}$ and $\hat{\upgamma} = \{ \upgamma'\geq \upgamma\}$
is the cone over $\upgamma$, 
then by directedness $\mathcal{A} =\{ \hat{\upgamma}\}$ satisfies the finite intersection
property and we will call 
\[ \mathfrak{c}=\mathfrak{c}_{\mathcal{I}}=\langle\mathcal{A}\rangle\] 
the {\bf {\small cone filter}} on $\mathcal{I}$. Note that $\mathfrak{c}$
is nonprincipal: indeed, if there were a set $X$ contained in all members of $\mathfrak{c}$, then for any $\upgamma_{0}\in X$
and $\upgamma>\upgamma_{0}$  we would have $X\not\subset\hat{\upgamma}\in \mathfrak{c}$.
An ultrafilter $\mathfrak{u}\supset\mathfrak{c}$ will
be called a  {\bf {\small cone ultrafilter}}\footnote{The ultraproduct proof of the compactness theorem of first order logic uses a cone ultrafilter on $\mathcal{I}={\sf Fin}(T)$ where $T$ is a finitely satisfiable first order
theory \cite{Po}.}.  Cone ultrafilters are nonprincipal.  Moreover, every element $X\in\mathfrak{u}$ of a cone ultrafilter is a directed set, and so
in particular, can be used to index nets.

The set of ultrafilters ${\sf Ult}(\mathcal{I})$ on $\mathcal{I}$, equipped with the topology 
generated
by the opens 
\[ V_{X}=\{ \mathfrak{u}\; |\;\; X\in \mathfrak{u} \} ,\quad X\in  {\sf 2}^{\mathcal{I}}   \]
is called the Stone space of $\mathcal{I}$ \cite{Jo}.  One has that $V^{\complement}_{X}= V_{X^{\complement}}$
where $\complement$ means complement, so that the $V_{X}$ are also closed.
With this topology, ${\sf Ult}(\mathcal{I})$ is totally-disconnected and compact, homeomorphic to 
the Stone-Cech compactification of $\mathcal{I}$ or dually, to 
the space of maximal ideals ${\sf Spec}({\sf 2}^{\mathcal{I}})$ equipped with the dual
Stone topology.  The isolated points are the principal ultrafilters.  

When $\mathcal{I}$ is directed, the subspace
${\sf Cone}(\mathcal{I})$ of cone ultrafilters is closed since
\[  {\sf Cone}(\mathcal{I}) = \bigcap_{\upgamma\in \mathcal{I}} V_{\hat{\upgamma}} . \]
In addition, ${\sf Cone}(\mathcal{I})$ is perfect as all of its elements are nonprincipal ultrafilters, hence are
non-isolated points.  In particular, ${\sf Cone}(\mathcal{I})$ is a (generalized) Cantor set, of cardinality possibly greater than that of the continuum.

\vspace{3mm}

\noindent  {\em Ultraproducts}

\vspace{3mm}

\noindent   
Let $L$ be a first order language, $\mathcal{I}$ a directed set and $\{ M_{\iota}\}_{\iota\in \mathcal{I}}$,
a family of $L$-structures
(e.g.\ a family of groups, rings, fields, {\it etc}) \cite{Ho}, \cite{Mark}.
Then the reduced product \cite{Ek} of 
the $M_{\iota}$ with respect to $\mathfrak{f}$ a filter on $\mathcal{I}$ is the
$L$-structure
\[ [M_{\iota}]_{\mathfrak{f}} := \prod M_{\iota}\big/\sim_{\mathfrak{f}} \]
where  $(x_{\iota} ) \sim_{\mathfrak{f}}( x_{\iota}')$ if and only if 
$   \{ \iota \; |\;\;  x_{\iota} = x_{\iota}' \}\in \mathfrak{f}$.  If
$M_{\iota}=M$ for all $\iota$, the reduced product is denoted
\[   \bast M_{\mathfrak{f}} \]
and called the reduced power of $M$ with respect to
$\mathfrak{f}$.
If $\mathfrak{f}=\mathfrak{u}$ is an ultrafilter, the reduced product
(reduced power) is called an ultraproduct (ultrapower).

By \L o\'{s}' Theorem \cite{Ho}, the ultrapower $\bast M_{\mathfrak{u}}$ is an elementary extension of $M$, where
the embedding $M\hookrightarrow \bast M_{\mathfrak{u}}$ is given by the constant nets. 
What this means is that $\bast M_{\mathfrak{u}}$ is a nonstandard model of $M$ i.e.\ it satisfies
the same set of first order $L$-sentences as $M$.  
In particular if $M$ is a group, ring or field than so is $\bast M_{\mathfrak{u}}$.
As one varies the ultrafilter, one obtains a sheaf 
\[ \bast\breve{M}\rightarrow {\sf Ult}(\mathcal{I}) \]
whose fiber over $\mathfrak{u}$ is $\bast M_{\mathfrak{u}}$, c.f. \cite{Mac}.

\begin{note}\label{CH} If one assumes the Continuum Hypothesis (CH) and $M$ is countable, then any two nonprincipal ultrafilters produce
isomorphic ultrapowers \cite{Ch-Kei}.    More generally, if the complete theory of $M$ is
uncountably categorical -- which is the case for $M=\C$ -- then again assuming the CH,
the nonprincipal fibers of $\bast\breve{M}$ will all be isomorphic, though not
canonically so \cite{Mark}.   We will not, however, assume CH
in this article.
\end{note}

If $\mathcal{I}=\N$ and $\mathfrak{u}$ is a fixed nonprincipal ultrafilter on $\N$, 
then we will suppress the ultrafilter in our notation and denote the ultrapower
\[ \bast M := \bast M_{\mathfrak{u}},\] informally referring to it as \guillemotleft nonstandard $M$\guillemotright ;
its elements will then be denoted $\bast x$, representatives of which are sequences $\{ x_{i}\}$ in $M$.

\vspace{3mm}

\noindent  {\em Extended reals}

\vspace{3mm}

\noindent We now turn to some specific ultrapowers which will be of interest to us:
the nonstandard versions of the integers,
the rationals, the reals and the complexes
related in the usual way: $\bast\Z\subset \bast\Q   \subset \bast\R \subset \bast\C $. 
Note that each of these structures contains classes corresponding to unbounded sequences,
and are therefore non-Archimedean (as rings or as fields). 
In addition, $\bast\Z$, $\bast\Q$ and $\bast\R$ are linearly ordered and the least upper bound
property does not hold in $\bast\R$ \cite{Ro}, \cite{Go}.

It can be easily checked that the field $\bast\Q$ is the field of fractions of the subring  $\bast\Z$.  In addition, $ \bast\Q$ is
also the field of fractions  of another, {\it local} subring, defined as follows.   Let $|\cdot  |$ be the Archimedean absolute value on $\Q$.
Then  $|\cdot  |$ induces in $\bast\Q$ a nonstandard absolute value with
values in $\bast\R_{+}$ = the nonnegative elements of $\bast\R$.  The set of bounded elements
\[  \bast\Q_{ {\rm fin}} = \{ \bast q\in\bast \Q\; |\;\; 
\text{there exists }r\in\R_{+} \text{ such that } |\bast q  |< r
  \}   \]
is a local ring with maximal ideal the set of infinitesimals
  \[  \bast\Q_{ \upvarepsilon} = \{ \bast q\in\bast \Q\; |\;\; 
\text{for all non-0 }r\in\R_{+}, \; |\bast q  |< r
  \} .  \]
We shall write $\bast x\simeq \bast y$ whenever
  $\bast x - \bast y\in \bast\Q_{ \upvarepsilon}$ and say that $\bast x$ and $\bast y$ are asymptotic or infinitesimal to one another.  We shall
  also refer to such a relation as an infinitesimal equation.
  
There is a canonical epimorphism 
  \[  {\rm std}: \bast\Q_{ {\rm fin}}\longrightarrow \R\]
  called the {\bf {\small standard part map}}: for any $\bast q\in \bast\Q_{ {\rm fin}}$,
  ${\rm std}(\bast q )$ is defined to be the unique accumulation point
  of any representative sequence $\{ q_{\upalpha}\}$ recognized by the ultrafilter.
More precisely, for any representative sequence $\{ q_{\upalpha}\}$, there exists $X\in\mathfrak{u}$
such that $\{ q_{\upalpha}\}|_{X}$ converges to a point ${\rm std} (\bast q )\in \R$, which depends 
neither on $\{ q_{\upalpha}\}$
nor on $X$ \cite{Ro}.

  The kernel of  ${\rm std}$ is $\bast\Q_{ \upvarepsilon}$ so that we have an isomorphism of fields
  \[   \bast\Q_{ {\rm fin}}/  \bast\Q_{ \upvarepsilon}  \cong \R . \]
  One may compare this situation with that of the $p$-adic numbers $\Q_{p}$,
  where the quotient of the ring of integers by its maximal ideal is the finite field $\F_{p}$ with $p$ elements.  

  Extending $|\cdot  |$ to $\bast\R$, we define in the same way the
  local ring $\bast \R_{\rm fin}$ with maximal ideal  $\bast \R_{\upvarepsilon}$ obtaining
  $\bast \R_{\rm fin}/\bast \R_{\upvarepsilon}\cong\R$.  We may similarly recover
  $\C$ from the quotient $\bast\C_{\rm fin}/\bast\C_{\upvarepsilon}$ 
  where $\bast\C_{\rm fin}$, $\bast\C_{\upvarepsilon}$ are defined using the usual absolute value in $\C$.

  The quotient
  \[  \bbull\R := \bast\R / \R_{\upvarepsilon} \cong \bast\Q /  \bast\Q_{ \upvarepsilon} \]
  is a real vector space (but not a topological vector space with respect to the quotient order topology) which
  we shall call the {\bf {\small extended reals}} \cite{Ge1},  \cite{Ge2}.  Note that $\bbull\R$ contains $\R$ canonically, and also $\bast\Z$ 
  since $\bast\Z\cap \bast \R_{\upvarepsilon} = \{0\}$.  We will view $\bbull\R$ as  \guillemotleft foliated\guillemotright\ by the
  cosets $\bbull x +\R$.  The subring $\bast\Z$ defines a transversal (in the sense that it has non trivial
  and discrete intersection with each coset leaf $\bbull x +\R$) and the  \guillemotleft leaf space\guillemotright\ may
be identified with $\bast\Z /\Z$, which a priori
is not endowed with any particular topology.  We define the extended complex numbers $\bbull\C$ in exactly the same way.

 \section{Diophantine Approximation Groups}\label{DAGs}

Fix $\bast\Z$ a nonstandard ring of integers and let $(\upmu , \uptheta)\in \pm\HP\times \SI^{1}$. 
Since $\Uplambda(\upmu )\subset\C$ is discrete, the ultrapower $\bast \Uplambda(\upmu )$
is naturally a subgroup of the vector space $\bbull \C=\bast \C/\bast \C_{\upvarepsilon}$.    In the Proposition which follows, we endow $\bbull\C $ with the the euclidean topology along its coset leaves $\bbull z +\C$
and the discrete topology transversally.
\begin{prop}\label{torusunif}
The quotient
\[ \bbull \C / \bast \Uplambda(\upmu )\]
is a topological group topologically isomorphic to $\T(\upmu )$.   
\end{prop}

\begin{proof} Note that  $(\bbull\C , +)$ is a 
topological group.  In addition,
$\bast\Uplambda(\upmu )$ is a complete transversal
for $\bbull\C$, so that every $\bbull z\in\bbull\C$ can be translated by an element of $\bast\Uplambda(\upmu )$
to $\C\subset\bbull\C$.  Since $\bast\Uplambda(\upmu )\cap \C =\Uplambda(\upmu )$, then
$\bbull \C / \bast \Uplambda(\upmu )= \C/\Uplambda(\upmu )$ and
 the Proposition follows.
\end{proof}

Define the {\bf {\small extended
line of $\boldsymbol\upmu$-slope $\boldsymbol\uptheta$}} as
\[ \bbull\widetilde{L}(\upmu ,\uptheta) :=\bbull\R\cdot ( \uptheta\upmu +1 )\subset \bbull \C . \]


\begin{defi}\label{dadef}  We say that $\bast n\in\bast\Z$ is a {\bf {\small diophantine
approximation}} of $\uptheta$ (relative to $\upmu$) if there exists $\bast m\in\bast\Z$
such that 
\begin{equation}\label{defcondition}
\bast m\upmu + \bast n\in \bbull\widetilde{L}(\upmu ,\uptheta)\cap \bast\Uplambda(\upmu )  .
\end{equation} 
\end{defi}

The next Proposition shows that the condition (\ref{defcondition}) depends only on $\uptheta$.

\begin{prop}\label{reldaisabsda}  The pair $(\bast m ,\bast n)$ defines a diophantine approximation of 
$\uptheta$ relative to $\upmu$ $\Leftrightarrow$ its coordinates satisfy  (in $\bbull\R$)
\begin{equation}\label{da}   \bast n \uptheta = \bast m 
 \end{equation} 
 \end{prop}
 
 \begin{proof}  The argument involves simple manipulations of equations {\it in} $\bbull\R$ using
 its $\R$-vector space structure. First, $\bast n$ is a diophantine
approximation of $\uptheta$ relative to $\upmu=a+ib$ $\Leftrightarrow$  there exists $\bbull r\in\bbull\R$ such that 
$ \bast m\upmu + \bast n = \bbull r (\uptheta\upmu +1)$ .  
Separating into real and imaginary parts gives 
\begin{equation}\label{pairda}
  \bast m a +\bast n = \bbull r (\uptheta a +1)\quad \text{and}\quad \bast m b=\bbull r\uptheta b .
   \end{equation}
The second equation of (\ref{pairda}) in turn yields $\bast m = \bbull r\uptheta $, which, when plugged back into the first
equation of (\ref{pairda}), gives
$   \bbull r \uptheta a +\bast n = \bbull r \uptheta a +\bbull r $
or $\bbull r = \bast n$.  Plugging the latter into $\bast m = \bbull r\uptheta $ gives (\ref{da}).  Conversely,
if $(\bast m ,\bast n)$ satisfies (\ref{da}), then taking $\bbull r = \bast n$ gives the pair of equations 
(\ref{pairda}), which imply the condition (\ref{defcondition}).
\end{proof}

The absolute version of diophantine approximation (\ref{da}) is that used in \cite{Ge1}, \cite{Ge3}.   It is clear from the form of (\ref{da}) that the collection of diophantine approximations of $\uptheta $ relative to $\upmu$ forms
a subgroup of $\bast\Z$ denoted 
\[ \bast\Z(\uptheta), \] which is independent of $\upmu$.   Note that this group is uncountably infinite and torsion-free.

\begin{theo} The group $\bast\Z(\uptheta)$ is an ideal in $\bast\Z$ $\Leftrightarrow$ $\uptheta\in\Q$.  If
$\uptheta, \uptheta '\in\SI^{1}$ satisfy $A(\uptheta )=\uptheta'$ for
some $A\in {\rm PGL}_{2}(\Z )$  then $\bast\Z(\uptheta )\cong\bast\Z(\uptheta ')$.
\end{theo}

This is proved in \cite{Ge1}; for the convenience of
the reader, we include here a 

\begin{proof}  If $\uptheta\in\Q$, then a pair $\bast n ,\bast m$ satisfies (\ref{da}) $\Leftrightarrow$
we have the equality {\it in the field} $\bast\R$: $ \bast n \uptheta = \bast m$.  Such an equality is invariant with respect
to multiplication by elements of $\bast\Z$, which shows that $\bast\Z(\uptheta)$ is an ideal.
If $\uptheta\in\R-\Q$ and 
$ \bast n\in\bast\Z(\uptheta ) $ then by irrationality we can find $\bast N\in \bast\Z$
such that $\bast N\bast n\uptheta$ contains a representative sequence asymptotic mod $\Z$ to any element of $\SI^{1}$ we choose.  If this element is not $0$, then 
$\bast N\bast n\not\in \bast\Z(\uptheta ) $, showing that  $\bast\Z(\uptheta)$ is not an ideal.
Now let \[ A=\left(\begin{array}{cc}
a & b \\
c & d \\
\end{array}
\right)\] be a (representative of an) element of  ${\rm PGL}_{2}(\Z )$.  If $(\bast m ,\bast n )$
satisfies (\ref{da}) then 
\[ \left(\begin{array}{c}
\bast m' \\
\bast n' 
\end{array}
\right) = A\left(\begin{array}{c}
\bast m \\
\bast n 
\end{array}
\right) = \left(\begin{array}{c}
 a \bast m+b\bast n \\
  c \bast m+d\bast n 
\end{array}
\right)  \]
satisfies the analogue of (\ref{da}) for $A(\uptheta )$.  Indeed, we have in $\bbull\R$ that 
$\bast n' A(\uptheta ) = \bast m'$ $\Leftrightarrow$
$\bast n' (a\uptheta +b) =\bast m' (c\uptheta +d) $ $ \Leftrightarrow$ $( c \bast m+d\bast n) (a\uptheta +b) =(a \bast m+b\bast n)(c\uptheta +d)$.  
But the latter equation is equivalent in $\bbull\R$
to $\bast n \uptheta =\bast m$.
\end{proof}

The element $\bast m$ associated
to $\bast n$ is unique: we refer to it as the dual of $\bast n$ and use the notation
\[ \bast n^{\perp} := \bast m.\]
The set of duals $\bast\Z^{\perp}(\uptheta)$ is a group, and when $\uptheta\not=0$, it is canonically isomorphic
to $\bast\Z(\uptheta)$ and equal to $\bast\Z(\uptheta^{-1})$.  In addition, 
\[ \bast\Uplambda(\upmu,\uptheta )= \bigg\{ \upmu\cdot\bast n^{\perp} +\bast n \bigg|\; \bast n\in \bast\Z(\uptheta)\bigg\}
\]
defines a subgroup of $\bast\Uplambda (\upmu )$
called the {\bf {\small group of $(\boldsymbol\upmu , \boldsymbol\uptheta )$-fractions}}, and the map $\bast n\mapsto  \upmu\cdot\bast n^{\perp} +\bast n$ defines
an isomorphism
$\bast\Z(\uptheta)\cong\bast\Uplambda(\upmu,\uptheta )$.  

The following Proposition expresses $ \bast\Uplambda(\upmu,\uptheta )$ as the intersection of a pair of ultrapowers of
standard lattices of rank 2 resp.\ 1.
Given $\upnu\in\pm\HP$, let 
\[  \bast \Uplambda_{\Updelta} (\upnu ) =\{ \bast n(1+\upnu )|\bast n \in\bast\Z \}\subset\bast\Uplambda(\upnu )\] 
be the ultrapower of the group $\Z\cdot(1+\upnu )\subset \Uplambda(\upnu )$ that uniformizes the diagonal cycle
$c(\upnu )\subset \T(\upnu )$.  
 
 \begin{prop}\label{latticeintersection} For any $\uptheta\in\R$, $\upmu\in\pm\HP$, 
\[   \bast\Uplambda(\upmu,\uptheta ) =    \bast\Uplambda(\upmu )\cap  \bast\Uplambda_{\Updelta} (\uptheta \upmu ).  \]
\end{prop}

\begin{proof}  By Proposition \ref{reldaisabsda}, if $\bast m\upmu+\bast n\in \bast\Uplambda(\upmu,\uptheta )$ then
$\bast m\upmu+\bast n=\bast n(\upmu \uptheta +1)\in  \bast\Uplambda_{\Updelta} (\uptheta \upmu )$.  Conversely
if $\bast m\upmu +\bast n= \bast n'(\upmu \uptheta +1)\in  \bast\Uplambda(\upmu )\cap\bast\Uplambda_{\Updelta} (\uptheta \upmu )$ then
$\bast m\upmu+\bast n\in \bbull\widetilde{L}(\upmu ,\uptheta)\cap \bast\Uplambda(\upmu )$.
\end{proof}

There is a natural homomorphism of abelian groups 
\[  {\rm std}(\upmu, \uptheta ): \bbull\widetilde{L}(\upmu ,\uptheta) \longrightarrow \T(\upmu ) \]
defined as follows. Take a representative sequence $\{ r_{\upalpha}\}\in \bbull r\in\bbull\R$,
and consider the image in $\T(\upmu )$ of the sequence 
\[ \{ r_{\upalpha}\cdot(\uptheta\upmu +1 )   \}\]
in the leaf $L(\upmu, \uptheta )\subset \T(\upmu )$ through the origin. Since $\T(\upmu )$ is compact, the ultrafilter
will recognize a unique limit point of this sequence, which is
independent of the choice of $\{ r_{\upalpha}\}\in \bbull r$ (again, see \cite{Ro} for more on this compactness principle).
We define 
\[ {\rm std}(\upmu, \uptheta )\big((\bbull r \cdot(\uptheta\upmu +1 ) \big)\] to be this limit point.

Notice that the leaves of $\bbull\R$ gives rise to leaves
of $\bbull\widetilde{L}(\upmu ,\uptheta)$, defined as the scalar multiples 
 $( \uptheta\upmu +1)\cdot (\bbull r+\R)$.  We note that the leaf corresponding to $\bbull r=0$ is the line $\widetilde{L}(\upmu ,\uptheta)\subset\C$
which was defined in \S 2: that is we have $\widetilde{L}(\upmu ,\uptheta)\subset\bbull\widetilde{L}(\upmu ,\uptheta)$.   The map ${\rm std}(\upmu, \uptheta )$ transports these leaves to the leaves of the associated Kronecker foliation $ \mathcal{F}(\upmu ,\uptheta)$.

\begin{theo}  If $\uptheta\in\R-\Q$ then ${\rm std}(\upmu, \uptheta )$ is surjective
with kernel $\bast\Uplambda(\upmu,\uptheta )$.
\end{theo}

\begin{proof}  Surjectivity follows from the density of $L(\upmu, \uptheta )$ in $\T(\upmu )$.
The map ${\rm std}(\upmu, \uptheta )$ coincides with the restriction of the epimorphism
$\bbull\C\rightarrow \T(\upmu )$ of Proposition \ref{torusunif} to the subspace $\bbull\widetilde{L}(\upmu ,\uptheta)$, 
so that the kernel is 
$\bbull\widetilde{L}(\upmu ,\uptheta) \cap  \bast\Uplambda(\upmu )= \bast\Uplambda(\upmu,\uptheta )$.
\end{proof}

Thus we have the \guillemotleft foliated group\guillemotright\ isomorphisms 
\[\bbull\R /\bast\Z(\uptheta )\cong \bbull\widetilde{L}(\upmu ,\uptheta)/\bast\Uplambda(\upmu,\uptheta )
 \cong \mathcal{F}(\upmu ,\uptheta).\]
Here we point out that the last isomorphism is not topological: nevertheless it is possible to put on $\bbull\widetilde{L}(\upmu ,\uptheta) $
a new transverse topology so that the action of $\bast\Uplambda(\upmu,\uptheta )$ is by homeomorphisms,
and that  the quotient $\bbull\widetilde{L}(\upmu ,\uptheta)/\bast\Uplambda(\upmu,\uptheta )$ becomes a foliation isomorphic
to $\mathcal{F}(\upmu ,\uptheta)$, see \cite{Ge1}, \cite{Ge3}.   By Proposition \ref{latticeintersection}, there are 
\guillemotleft covering maps\guillemotright 
\[ \mathcal{F}(\upmu ,\uptheta)\rightarrow \T(\upmu ) \quad \text{and}\quad
\mathcal{F}(\upmu ,\uptheta)\rightarrow c(\uptheta\upmu )\]
where $c(\uptheta\upmu )\subset \T(\uptheta\upmu )$ is the diagonal cycle.


\section{Ultrasolenoids}\label{ultrasolenoids}

As in previous sections we denote by $\bast\Z $, $\bast \C$ the ultrapowers
of $\Z$, $\C$ with respect to a fixed nonprincipal ultrafilter on $\N$.  
Let $S$ be a set, $\bast S_{\mathfrak{u}}$ be an arbitrary ultrapower with respect to some index
set $\mathcal{A}$ and ultrafilter $\mathfrak{u}$.  
Recall that a hyperfinite subset \cite{Go}, \cite{Ro} of $\bast S_{\mathfrak{u}}$ is an ultraproduct of the form
\[ [F_{\upalpha}]_{\mathfrak{u}}\subset  [S_{\upalpha} =S]_{\mathfrak{u}}  =\bast S_{\mathfrak{u}} \]
where $F_{\upalpha}\subset S$ is finite for all $\upalpha\in \mathcal{A}$.  Thus elements of $[F_{\upalpha}]$
are classes of sequences $\{x_{\upalpha}\}$ for which $x_{\upalpha}\in F_{\upalpha}$ for all $\upalpha$.
Note that every finite subset of $\bast S_{\mathfrak{u}} $ is hyperfinite.

Consider now the directed set (directed by inclusion)
\[\mathcal{H} = \big\{ \text{hyperfinite subsets } [F_{\upalpha}]\subset\bast\Z^{2}-\{0,0\}   \big\} .  \]
Denote by $\mathfrak{c}$ the cone filter on $\mathcal{H}$, which we recall was defined in \S 3 as 
the (nonprincipal) filter generated by the cones
\[ {\sf cone}\big([F_{\upalpha}]\big) = \bigg\{ [F'_{\upalpha}]\, \bigg|\; [F_{\upalpha}]\subset [F'_{\upalpha}] \bigg\} .\]
Let ${\sf Ult}(\mathcal{H})$ be the Stone space of ultrafilters on $\mathcal{H}$, and denote by ${\sf Cone}(\mathcal{H})\subset {\sf Ult}(\mathcal{H})$ the subspace of ultrafilters extending $\mathfrak{c}$.  
Each element $\mathfrak{u}\in {\sf Cone}(\mathcal{H})$ is nonprincipal, so ${\sf Cone}(\mathcal{H})$
is a Cantor set.  For us the importance of the cone ultrafilters is that they will provide partial summation
schemes that correspond well to the classical definition of a convergent infinite series.

Following \cite{Mac}, we define a sheaf $\bdiam\breve{\C}$ over ${\sf Ult}(\mathcal{H})$ as follows:
for each $\mathfrak{u}\in {\sf Ult}(\mathcal{H})$, the stalk over $\mathfrak{u}$, 
$\bdiam\C_{\mathfrak{u}}$, is the ultrapower
of $\bast\C$ with respect to $\mathfrak{u}$.
Let $\bdiam\breve{\Upgamma}$ be the $\bast\C$-algebra 
of set-theoretic sections of $\bdiam\breve{\C}$:  the $\bast\C$-algebra structure comes from the fact that $\bast\C$ is canonically included
in each fiber $\bdiam\C_{\mathfrak{u}}$ via the constant net inclusion.  In particular, we have canonical $\C$-algebra inclusions
$ \C \subset\bast \C \subset \bdiam\breve{\Upgamma}$
defined by the constant sections.
Finally, denote by $\bdiam\breve{\C}^{\rm cone}$
the restriction of $\bdiam\breve{\C}$ to ${\sf Cone}(\mathcal{H})$ and by
 $\bdiam\breve{\Upgamma}^{{\rm cone}}$
the sections of $\bdiam\breve{\C}^{\rm cone}$.  There is a canonical $\bast\C$-algebra epimorphism 
$\bdiam\breve{\Upgamma} \rightarrow \bdiam\breve{\Upgamma}^{{\rm cone}}$ given by restriction.

We now define subsheaves that correspond to 
$\uptheta\in\bar{\R}$.
For $\uptheta\not=\infty$ let
\[ \bast\Z^{2}(\uptheta) = \bigg\{(\bast n^{\perp}, \bast n) \bigg|\; \bast n\in \bast\Z(\uptheta)  \bigg\}
< \bast\Z^{2} .\]  
For $\uptheta=\infty$ we define $\bast\Z^{2}(\infty ):=\bast \Z\times\{ 0\}$.
Let $\mathcal{H}(\uptheta)\subset \mathcal{H}$
be the subset of hyperfinite subsets contained
in $\bast\Z^{2}(\uptheta)$. Let 
$\mathfrak{c}(\uptheta)$ be the
cone filter of $\mathcal{H}(\uptheta)$.
Denote by 
\[ {\sf Cone}(\mathcal{H})(\uptheta)\subset {\sf Ult}(\mathcal{H})\]
the subspace of ultrafilters $\mathfrak{u}$  of $\mathcal{H}$ ({\it not} of $\mathcal{H}(\uptheta)$) that extend $\mathfrak{c}(\uptheta)$.
The ultrafilters belonging to ${\sf Cone}(\mathcal{H})(\uptheta)$ are those ultrafilters of $\mathcal{H}$ that observe
the group $\bast\Z^{2}(\uptheta)$.

Let 
$\bdiam \breve{\C}^{\rm cone} (\uptheta)$ be
the restriction of $\bdiam \breve{\C}$ to  ${\sf Cone}(\mathcal{H})(\uptheta)$ and let
$\bdiam \breve{\Upgamma}^{\rm cone} (\uptheta) $ be
its $\bast\C$-algebra of sections.  The restriction map gives an algebra epimorphism
$\bdiam\breve{\Upgamma} \rightarrow \bdiam\breve{\Upgamma}^{{\rm cone}}(\uptheta )$. The Lemma which follows shows that the sheaves $\bdiam \breve{\C}^{\rm cone}$, $\bdiam \breve{\C}^{\rm cone} (\uptheta)$
are disjoint for all $\uptheta\in\bar{\R}$.  

\begin{lemm}\label{disjoint}  Let $\uptheta ,\eta\in \bar{\R}$ be distinct.  Then
\[ {\sf Cone}(\mathcal{H})(\uptheta)\cap {\sf Cone}(\mathcal{H}) = \emptyset   ={\sf Cone}(\mathcal{H})(\uptheta)\cap {\sf Cone}(\mathcal{H})(\upeta ) .\]
 \end{lemm}
\begin{proof} Suppose that $\mathfrak{u}\in {\sf Cone}(\mathcal{H})(\uptheta)\cap {\sf Cone}(\mathcal{H})$ so that
$\mathfrak{u}$ contains both $\mathfrak{c}$ and $\mathfrak{c}(\uptheta)$.  Let $F\subset\Z^{2}\subset
\bast\Z^{2}$ be a finite set containing a non zero element $(m,n)\not\in\bast\Z^{2} (\uptheta)$
(if $\uptheta\in\R-\Q$ this is true of any non-zero $(m,n)$).
Then ${\sf cone}(F)\in \mathfrak{c}\subset\mathfrak{u}$
and each element of ${\sf cone}(F)$ is a subset of $\bast\Z^{2}$ which contains  $(m,n)$.  
On the other hand, for any $X\in \mathfrak{c}(\uptheta)\subset\mathfrak{u}$, 
$X$ cannot contain any hyper finite subsets which contain $(m,n)$.
In particular we must have that $X\cap {\sf cone}(F)=	\emptyset$
is empty, contradicting the fact that $\mathfrak{u}$ is an ultrafilter.
Now for all $\uptheta\not=\eta$ we have 
$ \bast\Z^{2}(\uptheta)\cap \bast\Z^{2}(\eta ) = (0,0)$.
Indeed, if $(\bast n^{\perp}, \bast n)\in \bast\Z^{2}(\uptheta)$ then $\uptheta \simeq \bast n^{\perp}/ \bast n$
so that it is not possible to also have $\eta \simeq \bast n^{\perp}/ \bast n$ for $\uptheta\not=\eta$.
In particular $\mathcal{H}(\uptheta)\cap \mathcal{H}(\eta)=\{ (0,0)\}$ and therefore 
$ {\sf Cone}(\mathcal{H})(\uptheta)\cap {\sf Cone}(\mathcal{H})(\eta) =\emptyset$. 
\end{proof}


We define actions of ${\rm GL}_{2}(\Z )$ on the sheaves just considered, as well as on their algebras of sections.
First note that the left action of ${\rm GL}_{2}(\Z )$ on $\bast\Z^{2}-\{0,0\}$ induces one on hyperfinite sets, $[F_{\upalpha}]\mapsto [AF_{\upalpha}]$
for $A\in {\rm GL}_{2}(\Z )$.
This in turn induces 
an action on ${\sf Ult}(\mathcal{H})$ which preserves 
${\sf Cone}(\mathcal{H})$ and identifies ${\sf Cone}(\mathcal{H})(\uptheta)$ with ${\sf Cone}(\mathcal{H})(A(\uptheta ))$.

We can now define an action of ${\rm GL}_{2}(\Z )$ on
the sheaf $\bdiam \breve{\C}$ as follows: if $\big\{ \bast z_{[F_{\upalpha}]}\big\} $ represents an element
of $\bdiam z\in\bdiam\C_{\mathfrak{u}}$ then 
\[   \big\{ \bast w_{[F_{\upalpha}]}\big\}:=\big\{ \bast z_{[AF_{\upalpha}]} \big\} \]
represents an element of $\bdiam\C_{A^{-1}\mathfrak{u}}$.  
Indeed, suppose that $ \big\{\bast z'_{[F_{\upalpha}]}\big\} $ is another net representing $\bdiam z$.
Then there is a set $X\in\mathfrak{u}$ of hyperfinite sets such that 
$  \big\{\bast z_{[F_{\upalpha}]}\big\}\big|_{X} =  \big\{\bast z'_{[F_{\upalpha}]}\big\}\big|_{X}$.
It follows that
\[   \big\{ \bast w_{[F_{\upalpha}]}\big\}\big|_{A^{-1}X}  = 
 \big\{\bast z_{[F_{\upalpha}]}\big\}\big|_{X} =  \big\{\bast z'_{[F_{\upalpha}]}\big\}\big|_{X} =  \big\{ \bast w'_{[F_{\upalpha}]}\big\}\big|_{A^{-1}X}. \]
 Therefore $ \big\{ \bast w_{[F_{\upalpha}]}\big\}$ and $ \big\{ \bast w'_{[F_{\upalpha}]}\big\}$ are equivalent modulo $A^{-1}\mathfrak{u}$.
 
 Denote the action by $A\in{\rm GL}_{2}(\Z)$ defined in the previous paragraph by $A\odot \bdiam z$: we emphasize that it is 
{\it not} the matrix action along fibers.
Rather, the action is a shift, and so acts by $\bast\C$-algebra isomorphisms along the fibers of
$\bdiam \breve{\C}$ fixing the constant net classes.  That is, if $\bdiam z =\bast z\in \bdiam\C_{\mathfrak{u}}$ is a constant net class then
$A\odot \bast z = \bast z$ (viewed as an element of $ \bdiam\C_{A^{-1}\mathfrak{u}}$).
In particular, we obtain $\bast\C$-algebra isomorphisms $A:\bdiam\C_{A\mathfrak{u}}\rightarrow\bdiam\C_{\mathfrak{u}}$
for each $\mathfrak{u}\in{\sf Ult}(\mathcal{H})$.
This action stabilizes $\bdiam \breve{\C}^{\rm cone}$ and maps $\bdiam \breve{\C}^{\rm cone}(A(\uptheta ) )$
isomorphically onto $\bdiam \breve{\C}^{\rm cone}(\uptheta )$.

There is also an induced action by $A\in {\rm GL}_{2}(\Z)$ on elements of $\bdiam\breve{\Upgamma}$ defined
\[ f\mapsto g=A\odot f, \quad g(\mathfrak{u} ) := A\odot \big(f(A\mathfrak{u} )\big)  \]
which defines a $\bast\C$-isomorphism of $\bdiam \breve{\Upgamma}$ (since its acts as the identity on the
constant sections $\bast\C\subset \bdiam \breve{\Upgamma}$).  Again,  $\bdiam\breve{\Upgamma}^{{\rm cone}}$
is preserved by this action
and $ \bdiam\breve{\Upgamma}^{{\rm cone}}(A(\uptheta) )$ is identified with $ \bdiam\breve{\Upgamma}^{{\rm cone}}(\uptheta)$.

We now form the quotient of $\bdiam \breve{\C}$ with respect to the $ {\rm GL}_{2}(\Z)$ action:
\[ \bdiam \hat{\C} := {\rm GL}_{2}(\Z )\backslash \bdiam \breve{\C}.  \]
The result,
as such, is no longer a sheaf but rather a solenoid-like object, in a sense made precise in {\it Note} \ref{ultraidea} below.
We thus call $\bdiam \hat{\C}$ an {\bf {\small ultrasolenoid}}.  The quotient
\[ \bdiam \hat{\Upgamma}  := {\rm GL}_{2}(\Z )\backslash\bdiam \breve{\Upgamma} \]
is called the algebra of {\bf {\small ultratransversals}}: since ${\rm GL}_{2}(\Z)$ acts as the identity on the constant
sections, $ \bdiam \hat{\Upgamma}$ acquires the structure of a $\bast\C$-algebra extension of $\bast\C$.

In view of the fact that ${\rm GL}_{2}(\Z )$ stabilizes ${\sf Cone}(\mathcal{H})$, we also have a subultrasolenoid
\[ \hat{\C}^{\rm cl} :=  {\rm GL}_{2}(\Z )\backslash\bdiam \breve{\C}^{\rm cone}\subset\bdiam \hat{\C}\]
which we call the {\bf {\small classical ultrasolenoid}}
with the associated ultratransversal algebra:
\[  \bdiam \hat{\Upgamma}^{\rm cl}  := {\rm GL}_{2}(\Z )\backslash \bdiam \breve{\Upgamma}^{\rm cone} .\]
The restriction map $\bdiam\breve{\Upgamma}\rightarrow \bdiam\breve{\Upgamma}^{\rm cone}$ induces a projection
\[\uppi^{\rm cl}:  \bdiam \hat{\Upgamma}\longrightarrow \bdiam \hat{\Upgamma}^{\rm cl}.\]

Moreover, each equivalence class $[\uptheta]\in {\rm GL}_{2}(\Z)\backslash \bar{\R}$ gives rise as well to a subultrasolenoid
\[   \hat{\C}^{\rm qt}([\uptheta] ) :=  {\rm GL}_{2}(\Z )\big\backslash\bigg( \bigsqcup_{A\in{\rm GL}_{2}(\Z )}  \bdiam \breve{\C}^{\rm cone}(A\uptheta )\bigg)
\subset\bdiam \hat{\C} \]
which we call the {\bf [$\boldsymbol\uptheta$]}-{\bf {\small quantum ultrasolenoid}},
and a corresponding algebra of ultratransversals
\[   \hat{\Upgamma}^{\rm qt}([\uptheta] ) :=  {\rm GL}_{2}(\Z )\big\backslash\bigg( \bigsqcup_{A\in{\rm GL}_{2}(\Z )}  \bdiam \breve{\Upgamma}^{\rm cone}(A\uptheta )\bigg). \]
The restriction map induces again a projection
\[\uppi^{\rm qt}:  \bdiam \hat{\Upgamma}\longrightarrow \bdiam \hat{\Upgamma}^{\rm qt}([\uptheta]).\]
The motive for forming these ${\rm GL}_{2}(\Z)$ quotients is to ensure that the \guillemotleft Eisenstein objects\guillemotright\ we define in the sequel are automorphic.

\begin{note}\label{ultraidea}  The ultrasolenoids defined above are sheaf theoretic generalizations of the classical solenoid
$ \hat{\SI} = (\R\times\hat{\Z})/\Z$,
where $\hat{\Z}$ is the profinite completion of $\Z$ (a Cantor group) and the action is diagonal.
 We refer to the images of the fibers $\bdiam\C_{\mathfrak{u}}$ as the \guillemotleft leaves\guillemotright\ of
$\bdiam \hat{\C}$ (each of which is foliated by its $\bast\C$-cosets).   Recall from {\it Note} \ref{CH}
that assuming the Continuum Hypothesis,  the leaves of $\bdiam \hat{\C}$ will
be isomorphic to one another, though not canonically so.  This motivates the foliated view of these quotients.  In particular, the elements of
$\bdiam \hat{\Upgamma}$ are complete transversals of $\bdiam \hat{\C}$.
\end{note}
We describe briefly the context in which the constructions given in this section will be used
to define modular invariants.  As described in the above paragraphs we have a pair of epimorphisms
\begin{diagram}[heads=LaTeX]
& & \bdiam \hat{\Upgamma} & & \\
& \ldOnto^{\uppi^{\rm cl}} & &  \rdOnto^{\uppi^{\rm qt}([\uptheta] )} & \\
 \bdiam \hat{\Upgamma}^{\rm cl} & & & &\bdiam \hat{\Upgamma}^{\rm qt}([\uptheta])
\end{diagram}

We will define first in \S \ref{modinvariantsection} the universal modular invariant as a (set-theoretic) function
 \[ \bdiam\hat{\jmath}:\pm\HP\longrightarrow \bdiam \hat{\Upgamma}\]
for which 
\begin{itemize}
 \item[-]  the function $\bdiam\hat{\jmath}^{\rm cl}(\upmu ) := \uppi^{\rm cl}\circ\bdiam\hat{\jmath}(\upmu )$ yields a function asymptotic
 to the usual modular invariant of $\upmu $.
 \item[]
\item[-] the function $\bdiam\hat{\jmath}^{\rm qt}(\upmu ,\uptheta) := \uppi^{\rm qt}([\uptheta])\circ\bdiam\hat{\jmath}(\upmu)$ defines the (nonstandard) quantum
modular invariant of $(\upmu, \uptheta)$.  When $\upmu=i$, the result will be asymptotic to a multimap containing the standard quantum modular invariant defined in \S 
\ref{quantummodinvar}; a slight modification gives $j^{\rm qt}$ exactly.
\end{itemize}

We will then reinterpret the universal modular invariant as a {\it continuous} function 
 \[ \bdiam j :  \bdiam \widehat{\sf Mod}\longrightarrow \bdiam\hat{\C},\]
 where $\bdiam \widehat{\sf Mod}$ is a topological ultrasolenoid that we will define in the last part of 
 \S \ref{modinvariantsection}.

\section{Eisenstein Ultratransversals}\label{eisection}

We continue to fix as before ultrapowers $\bast\Z\subset \bast\Q\subset\bast\R\subset\bast\C$.
In this section we associate to every $k\in\Z$ and each pair $(\upmu , \uptheta)\in\pm\HP\times \SI^{1}$
an analogue of the classical Eisenstein series, defined 
as an ultratransversal
 \[  \bdiam\hat{G}_{k}^{\rm qt}(\upmu ,\uptheta )\in \bdiam \hat{\Upgamma}^{\rm qt} .\]

Fix $\upmu\in\pm\HP$.  For each hyperfinite set
$[F_{\upalpha}]\in \mathcal{H}$ and $k\in\Z$ consider the hyperfinite sum \cite{Go}
\[ G_{k }(\upmu )_{[F_{\upalpha}]} = \sum_{[F_{\upalpha}]} ( \bast m\upmu +\bast n)^{-2k} :=
\text{ $\ast$-class of } \left\{ \sum_{(m_{\upalpha}, n_{\upalpha})\in F_{\upalpha}} ( m_{\upalpha}\upmu + n_{\upalpha})^{-2k}\right\}
 \;\in\bast\C . \]
 Note that this expression is well-defined even for $k<2$, in contrast with to the classical situation.  For example, when
 $k=0$, we have 
 \[ G_{0 }(\upmu )_{[F_{\upalpha}]} = \text{ hypercardinality of } [F_{\upalpha}] =  \ast\text{-class of }\{ |F_{\upalpha} | \}.\]

The $\mathcal{H}$-net 
\begin{equation}\label{definingnet} \left\{ G_{k }(\upmu )_{[F_{\upalpha}]}  \right\} 
\end{equation}
defines as described in the previous section an element
$ \bdiam\hat{G}_{k }(\upmu )\in \bdiam\breve{\Upgamma}. $
We thus obtain a function
\[ \bdiam\hat{G}_{k}:\pm \HP\longrightarrow \bdiam\hat{\Upgamma} .\]
\begin{prop}\label{automorphyofthecl} $\bdiam\hat{G}_{k}$ is a modular form of weight $k$:
\[ (A'(\upmu))^{k}\cdot \bdiam\hat{G}_{k}(A\upmu) =  \bdiam\hat{G}_{k}(\upmu ) \]
for all $A\in {\rm GL}_{2}(\Z )$.
\end{prop}

\begin{proof}  We will show that in $\bdiam\breve{\Upgamma}$
\[  (A'(\upmu))^{k} \cdot \bdiam\breve{G}_{k}(A\upmu) =   A^{T} \odot \bdiam\breve{G}_{k }(\upmu ) . \]
We calculate at the level of  the net (\ref{definingnet}): for
\[ A = \left( \begin{array}{ll} 
                     a & b \\
                     c & d
              \end{array}
  \right)\in {\rm GL}_{2}(\Z )  \]
we have
\begin{eqnarray*}
(A'(\upmu))^{k} \cdot G_{k }(A\upmu )_{[F_{\upalpha}]} & = & (c\upmu + d)^{-2k}\sum_{[F_{\upalpha}]}\bigg(\bast m\bigg(\frac{a\upmu +b}{c\upmu +d} \bigg)+
\bast n\bigg)^{-2k} \\ 
& \\
& = &  \sum_{[F_{\upalpha}]}\big((a\bast m + c\bast n)\upmu +(b\bast m+d\bast n) \big)^{-2k} \\
& \\
& = & G_{k }(\upmu )_{A^{T}[F_{\upalpha}]}
\end{eqnarray*}
from which the statement follows.
\end{proof}

Denote by 
\[ \bdiam\hat{G}^{{\rm cl}} _{k}: \pm\HP\longrightarrow  \bdiam \hat{\Upgamma}^{\rm cl}  \]
the composition of $\bdiam\hat{G}_{k}$ with the projection
$\uppi^{\rm cl} :\bdiam\hat{\Upgamma}\longrightarrow  \bdiam \hat{\Upgamma}^{\rm cl} $ defined in \S \ref{ultrasolenoids}.
For $k\geq 2$ and $\upmu\in\pm\HP$, let $G_{k }(\upmu )$ be the usual (standard) Eisenstein series.  Since $\C\subset \bdiam \hat{\Upgamma}^{\rm cl}$
we may view $G_{k }$
as defining a family of \guillemotleft constant
ultratransversals\guillemotright
\[ G_{k}:\HP\rightarrow \bdiam \hat{\Upgamma}^{\rm cl}. \] 

In what follows, for any pair of sections $f,g:{\sf Ult}(\mathcal{H})\rightarrow\bdiam\breve{\C}$ we write 
$f\simeq g$ if $f-g\in\bast\C_{\upvarepsilon}\subset\bast\C\subset \bdiam\C_{\mathfrak{u}}$ for all
$\mathfrak{u}$.  Notice that this relation is preserved by the action of ${\rm GL}_{2}(\Z)$,
giving rise to the relation of infinitesimality of ultratransversals in $\bdiam \hat{\Upgamma}^{\rm cl}$.

\begin{prop}\label{asymtoclass}  For all $\upmu\in\pm\HP$, 
\[ G_{k}(\upmu)\;\simeq\; \bdiam\hat{G}_{k }(\upmu )^{\rm cl} .\]
\end{prop}

\begin{proof}  Let $\mathfrak{u}\in {\sf Cone}(\mathcal{H})$.  It will be enough to check that
for any finite subset $F\subset\Z^{2}$ that the net of hyperfinite sums over elements
in ${\sf cone}(F)$ converges to $G_{k }(\upmu )$.  This is certainly true if we restrict
to the subnet of all {\it finite} subsets $F'\supset F$, because the classical Eisenstein series
converges.  Now if $[F_{\upalpha}]\supset F$ is a general hyperfinite containing $F$, and $\upvarepsilon>0$, let $F'\supset F$
be such that $G_{k }(\upmu )_{F''}$ is $\upvarepsilon$-close to $G_{k}(\upmu)$ for all $F''\supset F'$. Define
$F_{\upalpha}'=F_{\upalpha}\cup F'$.  Then $[F_{\upalpha}']\supset [F_{\upalpha}]$ and $G_{k }(\upmu )_{[F_{\upalpha}']}$ has standard part which is $\upvarepsilon$-close to $G_{k}(\upmu)$; moreover, every $[F_{\upalpha}'']\supset [F_{\upalpha}']$ has the same property.  It follows that the net of hyperfinite sums
associated to ${\sf cone}(F)$
have standard parts converging to $G_{k}(\upmu)$.  The infinitesimality statement follows.
\end{proof}

The proof of the Proposition \ref{asymtoclass} reveals the function of cone ultrafilters: they are the
ones that recognize classically convergent infinite series.

Define
\[ \bdiam\hat{G}_{k}^{\rm qt}:\pm\HP\times\SI^{1}\longrightarrow  \bdiam \hat{\Upgamma}^{\rm qt},
\quad (\upmu ,\uptheta )\mapsto \uppi^{\rm qt}([\uptheta]) \big(  \bdiam\hat{G}_{k }(\upmu  ) \big).\]
With the action of ${\rm GL}_{2}(\Z )$ on $\pm\HP\times \SI^{1}$ defined
as in (\ref{Aactiononpairs}) we have 

\begin{prop} $\bdiam\hat{G}_{k}^{\rm qt}$ is a modular form of weight $k$:
\[ (A'(\upmu))^{k}\cdot \bdiam\hat{G}_{k}^{\rm qt}(A(\upmu ,\uptheta )) =  \bdiam\hat{G}^{\rm qt}_{k}(\upmu ,\uptheta ).\]
\end{prop}

\begin{proof}  Exactly the same proof as Proposition \ref{automorphyofthecl}.
\end{proof}

Let $\bdiam\breve{\R}\subset \bdiam\breve{\C}$ be the sheaf of real points, and denote by
$\bdiam\breve{\Upgamma}(\R )$ the sections with values in $\bdiam\breve{\R}$.  
Let $ \bdiam \hat{\Upgamma}^{\rm qt}(\R )$ denote the associated real points in $ \bdiam \hat{\Upgamma}^{\rm qt}$.  
For the value $\upmu =i$, it is well-known that the classical Eisenstein series is real
valued.  For the same reasons we have the following
important reality result for the classical Kronecker foliations (those corresponding to
pairs $(i,\uptheta )$):

\begin{prop}\label{realocus}  For all $k$ and $\uptheta\in\bar{\R}$,
$ \bdiam\hat{G}_{k}^{\rm qt}(i,\uptheta )\in  \bdiam \hat{\Upgamma}^{\rm qt}(\R )$.
\end{prop}

\begin{proof}  As before, we work on the level of the defining net (\ref{definingnet}). We consider any subnet
\[ \bigg\{ \sum_{[F_{\upalpha}]} (\bast n +\bast n^{\perp}i)^{-2k} \bigg\} \]
where $[F_{\upalpha}]$ range over the elements of some $X\in\mathfrak{c}(\uptheta)$.  Taking the
conjugate yields 
\[  \big\{  \sum_{[F_{\upalpha}]} (\bast n -\bast n^{\perp}i)^{-2k}  \big\} =  \big\{ \sum_{[F_{\upalpha}]} (\bast n +\bast n^{\perp}A(i))^{-2k}  \big\} \]
where $A$ is the element of ${\rm PGL}_{2}(\Z )$ defining $z\mapsto -z$.  It follows then by
the automorphy that 
$ \bdiam\hat{G}_{k}^{\rm qt}(i,\uptheta )$ is equal to its own conjugate. 
\end{proof}

\section{The Universal Modular Invariant}\label{modinvariantsection}

In this section we will define a universal modular invariant as a map of ultrasolenoids, in such a way that each of the classical
and quantum invariants may be recovered from it as a subquotient (a restriction followed by quotients).

Define the {\bf {\small universal modular invariant}} 
\[  \bdiam \hat{\jmath}:\pm\HP\rightarrow \bdiam \hat{\Upgamma}\] via the classical template \cite{Hu}:
\[ \bdiam \hat{\jmath}(\upmu)= 12^{3}\cdot \frac{   \bdiam \hat{g}_{2}(\upmu)^{3} }{ \bdiam \hat{g}_{2}(\upmu )^{3} - 27 \cdot \bdiam \hat{g}_{3}(\upmu )^{2}},  \]
where the lower case (normalized) Eisenstein ultratransversals $ \bdiam \hat{g}_{2}$, $ \bdiam \hat{g}_{3}$ are defined in the usual way by scaling  $\bdiam \hat{G}_{2},  \bdiam \hat{G}_{3}$ by 60 resp.\ 140.


The classical and quantum modular
invariants are defined by composition with the projections $\uppi^{\rm cl}$ resp.\ $\uppi^{\rm qt} (\uptheta )$ (see the end of \S \ref{ultrasolenoids}), making the quantum invariant
a function of $\pm\HP\times \SI^{1}$;  since the automorphies of the numerator and denominator cancel, we obtain modular functions
 \[  \bdiam  \hat{\jmath}^{\rm cl} :=\uppi^{\rm cl}\circ \bdiam \hat{\jmath} : {\sf Mod}^{\rm cl}\longrightarrow \bdiam \hat{\Upgamma}^{\rm cl} \]
 and 
 \[  \bdiam  \hat{\jmath}^{\rm qt} : {\sf Mod}^{\rm kf} \longrightarrow \bdiam \hat{\Upgamma}^{\rm qt}, \quad
  \bdiam  \hat{\jmath}^{\rm qt} (\upmu,\uptheta)= \uppi^{\rm qt}([\uptheta ])\circ \bdiam \hat{\jmath}(\upmu )
  .\]
By Proposition \ref{asymtoclass}, $\bdiam \hat{g}^{\rm cl}_{2}\simeq g_{2}$ and $\bdiam \hat{g}^{\rm cl}_{3}\simeq g_{3}$
so we
have immediately:

\begin{coro}\label{classasympclass}  Let $j(\upmu )$ be the usual modular invariant of the elliptic curve $\T(\upmu )$
viewed as a constant transversal in $ \bdiam \hat{\Upgamma}^{\rm cl} $.  Then $\bdiam  \hat{\jmath}^{\rm cl} (\upmu )\simeq j(\upmu )$. 
\end{coro}

Note that by Proposition \ref{realocus}, the image of  ${\sf Mod}^{\rm kf}_{[i]}$ (the 
fiber over $[i]$) by $ \bdiam  \hat{\jmath}^{\rm qt}$ belongs to the real locus:

\begin{prop}\label{jreal} For all $\uptheta\in\bar{\R}$, $ {}^{\diamond}\hat{\jmath}^{\rm qt}(i,\uptheta )\in {}^{\diamond}\hat{\Upgamma}^{\rm qt}(\R )$.
\end{prop}

The reality of the image of  ${\sf Mod}^{\rm kf}_{[i]}$ given by Proposition \ref{jreal} 
suggests that we may calculate $ {}^{\diamond}\hat{\jmath}^{\rm qt}(i,\uptheta )$ using hyperfinite partial sums over the group $\bast\Z(\uptheta)\subset\bast\R$
rather than over the group $\bast\Uplambda(i,\uptheta)\subset\bast\C$. 
This is reasonable, since by Proposition \ref{reldaisabsda}, every element of $\bast\Uplambda (i, \uptheta )$ is of the form
 $\bast n^{\perp}i+\bast n$ for $\bast n\in\bast\Z(\uptheta)$.

We recall that the ultratransversal $ {}^{\diamond}\hat{\jmath}^{\rm qt}(i,\uptheta)$ is an equivalence class
of the section $ {}^{\diamond}\breve{\jmath}^{\rm qt}(i,\uptheta)$ of the sheaf ${}^{\diamond}\breve{\C}^{\rm qt}$, where for each $\mathfrak{u}\in{\sf Cone}(\mathcal{H})(\uptheta )$,
the value of $ {}^{\diamond}\breve{\jmath}^{\rm qt}(i,\uptheta)$ is the
$\mathfrak{u}$-class of the net of hyper-finite partial sums 
$    \big\{    j_{[F_{\upalpha}] }\big\}$, for $[F_{\upalpha}] \subset \bast\Uplambda (i, \uptheta )$  hyper-finite. 
Let us write
 \[   j_{F_{\upalpha}} (i,\uptheta ) := \frac{12^{3}}{ 1-J_{F_{\upalpha}}(i,\uptheta ) }  ,\quad J_{F_{\upalpha}}(i,\uptheta ) = \frac{49}{20}\frac{\big( \sum_{mi+n\in F_{\upalpha}}  (mi+n)^{-6} \big)^{2}   }{\big( \sum_{mi+n\in F_{\upalpha}}  (mi+n)^{-4} \big)^{3}} . \]
 Since we are considering nets of hyperfinite partial sums which are increasing w.r.t.\ inclusion i.e.\ 
nets indexed by $X\in\mathfrak{c}(\uptheta)$,
we may assume that the hyperfinite set $[F_{\upalpha}]$ is symmetric w.r.t.\ multiplication by $-1$, and write
 \[  J_{F_{\upalpha}}(i,\uptheta ) =  \frac{49}{40} \frac{\big( \sum_{mi+n\in F_{\upalpha},\; n>0}  (mi+n)^{-6} \big)^{2}   }{\big( \sum_{mi+n\in F_{\upalpha},\; n>0}  (mi+n)^{-4} \big)^{3}} .\]
 For such a hyperfinite set $[F_{\upalpha}]$, let $[\bar{F}_{\upalpha} ]\subset \bast\Z(\uptheta )$ be the set of $\bast n\in \bast\Z(\uptheta )$
 for which $\bast n^{\perp}i + \bast n\in [F_{\upalpha}]$.  Consider the hyperfinite sum
 $j_{[\bar{F}_{\upalpha}]}(\uptheta )$ defined as the $\bast$-class corresponding to the sequence
 \[  j_{\bar{F}_{\upalpha}} (\uptheta ) := \frac{12^{3}}{ 1-J_{\bar{F}_{\upalpha}}(\uptheta ) }, \quad
J_{\bar{F}_{\upalpha}}(\uptheta ) =  \frac{49}{40} \frac{\big( \sum_{n\in \bar{F}_{\upalpha},\; n>0}  n^{-6} \big)^{2}   }{\big( \sum_{n\in \bar{F}_{\upalpha},\; n>0}  n^{-4} \big)^{3}} . 
      \]
 
 \begin{prop}\label{ABProp}  There exists $\bast u\in\bast \R$, $\bast u\simeq 1$, with \[ j_{[F_{\upalpha}]}(i,\uptheta )=\bast u\cdot  j_{[\bar{F}_{\upalpha}]}(\uptheta ).\]  In particular,
 $j_{[F_{\upalpha}]}(i,\uptheta )\in \bast\R_{\rm fin}$ $\Leftrightarrow$ $j_{[\bar{F}_{\upalpha}]}(\uptheta )\in \bast\R_{\rm fin}$ 
 and in this case,  $j_{[F_{\upalpha}]}(i,\uptheta )\simeq j_{[\bar{F}_{\upalpha}]}(\uptheta )$.
 \end{prop}
 
 \begin{proof} Let $[F_{\upalpha}]\subset  \bast\Uplambda (i, \uptheta )$; then we may write for each $\upalpha$
 \[  J_{F_{\upalpha}}(i,\uptheta) = \frac{49}{40} \frac{\bigg( \sum_{mi+n\in F_{\upalpha},\; n>0}  n^{-6}\big((m/n)i+1\big)^{-6} \bigg)^{2}   }{\bigg( \sum_{mi+n\in F_{\upalpha},\; n>0}  n^{-4}\big((m/n)i+1\big)^{-4} \bigg)^{3}}  .\]
 Since $[F_{\upalpha}]\subset  \bast\Uplambda (i, \uptheta )$, for any $\upvarepsilon >0$ we have $|\uptheta - m/n|=\upvarepsilon(n)<\upvarepsilon$ for all $mi+n\in F_{\upalpha}$ and $\upalpha$ sufficiently large.  Multiplying
 numerator and denominator by $(\uptheta i+1)^{12}$ and writing 
 \[ (m/n)i +1= \uptheta i +1 \pm \upvarepsilon (n)\] gives the bounds
 \[ \left( \frac{\uptheta + 1 }{\uptheta + 1 +\upvarepsilon}\right) ^{12} J_{\bar{F}_{\upalpha}}(\uptheta ) < J_{F_{\upalpha}}(i,\uptheta) < \left( \frac{\uptheta + 1}{\uptheta + 1 -\upvarepsilon}\right) ^{12} J_{\bar{F}_{\upalpha}}(\uptheta ) .\]  
 The result follows immediately.
 \end{proof}

 
 Let $\bdiam\hat{\sigma}\in  \bdiam\hat{\Upgamma}^{\rm qt}(\R )$ be any section class.  We say that  $\bdiam\hat{\sigma}$
 has standard part at (the ${\rm GL}_{2}(\Z )$ orbit of) $\mathfrak{u}$
 \[ {\rm std}\big( \bdiam\hat{\sigma}(\mathfrak{u})\big)\in\R\]
 if there exists a representative section $\bdiam\breve{\sigma}$ such that
 for all $M\in\R_{+}$, 
 \[  \big|  \bdiam\breve{\sigma}(\mathfrak{u}) - {\rm std}\big( \bdiam\hat{\sigma}(\mathfrak{u})\big)    \big| <M. \]
Notice that if the standard part at $\mathfrak{u}$ exists it is unique.
 If $\bdiam\hat{\sigma}$
does not have standard part at (the ${\rm GL}_{2}(\Z )$ orbit of) $\mathfrak{u}$
we will write 
\[{\rm std}\big( \bdiam\hat{\sigma}(\mathfrak{u})\big)=\infty.\]
Thus each section class and each $\uptheta\in\R$ determines, in particular, a function on the ${\rm GL}_{2}(\Z )$-orbit of
${\sf Cone}(\mathcal{H})(\uptheta )$
\[ {\rm std}(\bdiam\hat{\upsigma})_{\uptheta}:
 {\rm GL}_{2}(\Z )\bigg\backslash \bigg(\bigsqcup_{A\in {\rm GL}_{2}(\Z)} {\sf Cone}(\mathcal{H})(A(\uptheta ))\bigg)\longrightarrow \bar{\R} .\]
This gives an induced multimap
\[  \overline{\rm std}(\bdiam\hat{\upsigma}): \R\multimap \bar{\R}, \quad
  \uptheta\mapsto {\rm std}(\hat{\upsigma})_{\uptheta} \big({\sf Cone}(\mathcal{H})(\uptheta )\big).\]
 
 
 \begin{theo}\label{stdprtlim} Let $\uptheta\in\R-\Q$.  Then $\overline{\rm std}({}^{\diamond}\hat{\jmath}^{\rm qt}(i,\uptheta))\supset j^{\rm qt}(\uptheta )$.
 \end{theo}
 
 \begin{proof}  Suppose that $ j_{0}\in j^{\rm qt}(\uptheta ) $ is the limit corresponding to the sequence $\{ \upvarepsilon_{\upalpha}\}$, whose class in
 $\bast\R_{\upvarepsilon}$ is denoted $\bast\upvarepsilon$.  Consider 
 a {\it shift function} $\upsigma: \N\rightarrow\N$ i.e.\ a function which is finite-to-1 and does not reverse the order: $\upsigma (\upalpha )\leq \upsigma (\upbeta )$
 if $\upalpha\leq\upbeta$.  For such a shift function, the sequence $  \{ \upvarepsilon_{\upsigma (\upalpha)}\}$ will produce the same 
  limit $j_{0}$.  We denote the class of such a shifted sequence by $\upsigma (\bast\upvarepsilon )$.\footnote{Note that the map $\upsigma$ {\it does not} induce a function on $\bast\R$.}  Note that if $\bast\updelta$ is any positive infinitesimal  
 there exist shifts $\upsigma_{0}, \upsigma_{1}$ with $\upsigma_{0} (\bast\upvarepsilon ) <\bast\updelta<\upsigma_{1} (\bast\upvarepsilon )$.
 Let 
 \[  \bast\Z^{2}_{\bast \upvarepsilon}(\uptheta ):=\big\{ (\bast m,\bast n)\in\bast\Z^{2}(\uptheta )|\; \pm\bast n \in [ B_{\upvarepsilon_{\upalpha}}(\uptheta )]
  \big\} \] 
  where $[ B_{\upvarepsilon_{\upalpha}}(\uptheta )]$ is the ultraproduct of the $B_{\upvarepsilon_{\upalpha}}(\uptheta )\subset\N$ defined in \S \ref{quantummodinvar}.
 We will produce a set of hyperfinites $\sf{X}\subset\mathcal{H}$ compatible with the cone filter $\mathfrak{c}(\uptheta )$ (i.e.
 ${\sf X}\cap {\sf Cone}([F_{\upalpha}])\not=\emptyset$ for all $[F_{\upalpha}]\subset \bast \Uplambda (i,\uptheta )$) for which the net
 $\{ j_{[F_{\upalpha}]}\}_{[F_{\upalpha}]\in{\sf X}}$ converges to $j_{0}$.  Given $[F_{\upalpha}]\subset \bast \Z^{2} (\uptheta )$ we can find some
 shift map $\upsigma$ so that $[F_{\upalpha}]\subset\bast\Z^{2}_{\upsigma(\bast \upvarepsilon)}(\uptheta )$.  We may choose a hyperfinite
 $[\tilde{F}_{\upalpha}]\subset \bast\Z^{2}_{\upsigma(\bast \upvarepsilon)}(\uptheta )\cap {\sf Cone}([F_{\upalpha}])$ such that
 $j_{[\tilde{F}_{\upalpha}]}\simeq j_{0}$.  Let ${\sf X}$ be the set of such $[\tilde{F}_{\upalpha}]$.  Then ${\sf X}$ is compatible with 
 $\mathfrak{c}(\uptheta )$ so there exists a cone ultrafilter $\mathfrak{u}$ containing ${\sf X}$.  This ultrafilter will produce a standard part $j_{0}$.

\end{proof}

In order to recover the invariant on the nose requires a slight paring down of ${\sf Cone}(\uptheta )$.  Consider
the subset ${\sf J}(\uptheta )\subset \mathcal{H}(\uptheta )$ of hyperfinites $[F_{\upalpha}]\subset \bast\Z^{2}(\uptheta )$ for which
\[  j_{[F_{\upalpha}]}(i,\uptheta )\simeq j_{0}\in j^{\rm qt}(\uptheta ).  \]
Clearly ${\sf J}(\uptheta )$ intersects nontrivially every element of $\mathfrak{c}(\uptheta )$ (it intersects all the cones contained in
 $\mathfrak{c}(\uptheta )$), therefore we may form
the filter $\mathfrak{c}(\uptheta )_{0}$ generated by $\mathfrak{c}(\uptheta )$ and ${\sf J}(\uptheta )$.  The set of
ultrafilters $\mathfrak{u}$ on $\mathcal{H}$ extending $\mathfrak{c}(\uptheta )_{0}$ defines a closed subset
\[ {\sf Cone}(\uptheta )_{0}\subset {\sf Cone}(\uptheta ).\]
Note that ${\sf Cone}(\uptheta )_{0}$ is taken to ${\sf Cone}(A(\uptheta) )_{0}$ by any $A\in{\rm GL}_{2}(\Z )$,
since $A({\sf J}(\uptheta ))={\sf J}(A(\uptheta) )$: the latter follows from the fact that the standard $j^{\rm qt}$
is ${\rm GL}_{2}(\Z )$-invariant. 
 If we denote by ${}^{\diamond}\hat{\jmath}^{\rm qt}_{0}(i,\uptheta)$ the restriction of ${}^{\diamond}\hat{\jmath}^{\rm qt}(i,\uptheta)$
 to $ {\sf Cone}(\uptheta )_{0}$ then we have
 
 \begin{theo}\label{stdprtlim} Let $\uptheta\in\R-\Q$.  Then $\overline{\rm std}({}^{\diamond}\hat{\jmath}^{\rm qt}_{0}(i,\uptheta))= j^{\rm qt}(\uptheta )$.
 \end{theo}
 
 \begin{proof} By the previous Theorem we only need to show that every point in $\overline{\rm std}({}^{\diamond}\hat{\jmath}^{\rm qt}_{0}(i,\uptheta))$ belongs to $j^{\rm qt}(\uptheta )$.  If $j_{0}$ is a standard part of ${}^{\diamond}\hat{\jmath}^{\rm qt}_{0}(i,\uptheta)$ 
converging with respect to the restriction to ${\sf X}\in\mathfrak{u}$, then we may produce
an element $[F_{\alpha}]\in{\sf J}(\uptheta )\cap {\sf X}$ with $ j_{[F_{\upalpha}]}(i,\uptheta )\simeq j_{0}$ and by definition of ${\sf J}(\uptheta )$, this shows that $j_{0}\in  j^{\rm qt}(\uptheta )$.
 \end{proof}

\vspace{3mm}

\begin{center}
$\maltese$
\end{center}

\vspace{3mm}

We finish with an equivalent formulation of $\bdiam\hat{\jmath}$ that evokes the framework of the classical  modular invariant.  
Let 
\[ \pm\bdiam\breve{\HP}\subset \bdiam\breve{\C}\] be the subsheaf of $ \bdiam\breve{\C}$ having stalk $\pm\bdiam\HP_{\mathfrak{u}}\subset\bdiam\C_{\mathfrak{u}}$ for each $\mathfrak{u}\in \text{\sf Ult}(\mathcal{H})$.
On $\pm\bdiam\breve{\HP}$ we define a diagonal action of ${\rm GL}_{2}(\Z )$, given, at the level of nets, by
\[ \left\{\bast z _{[F_{\upalpha}]}\right\} \mapsto \left\{\bast w _{[F_{\upalpha}]}\right\}:=\left\{A\left(\bast z_{[A^{T}F_{\upalpha}]}\right)\right\}.\]
The action of $A$ is therefore a shifted linear action.
We denote it by 
\[  \bdiam z_{\mathfrak{u}} \mapsto A\circledast \bdiam  z_{\mathfrak{u}}\in \C_{A^{-T}\mathfrak{u}}\]
to distinguish it from the earlier defined $A\odot \bdiam z_{\mathfrak{u}}$, the shift induced by $A$.
The quotient by this action is an ultrasolenoid denoted
\[  \bdiam \widehat{\sf Mod} , \]
whose \guillemotleft leaves\guillemotright\ are the images of the stalks of $\pm\bdiam\breve{\HP}$.

The action which defines $ \bdiam \widehat{\sf Mod} $ is the analogue of the diagonal action of ${\rm GL}_{2}(\Z )$ on $\pm\HP\times\SI^{1}$ used to
define the signed Anosov foliation ${\sf Mod}^{\rm kf}$.
We may think of each point of $\bdiam \widehat{\sf Mod}$ as parametrizing an isomorphism class of 
abstract nonstandard Kronecker foliation
in the nonstandard elliptic curve 
\[ \bdiam\C_{\mathfrak{u}}\,/\;\bdiam \Uplambda(\bdiam z_{\mathfrak{u}}),\] where
$\bdiam \Uplambda(\bdiam z_{\mathfrak{u}})$ is the $\bdiam\Z_{\mathfrak{u}}$-module generated by $1,\bdiam z_{\mathfrak{u}}$
 and where moreover
each $\mathfrak{u}\in \text{\sf Ult}(\mathcal{H})$ is to be thought of as supplying the data of an  \guillemotleft ultraslope\guillemotright.  

Denote by $\bdiam{\sf Mod}$ the image in $\bdiam \widehat{\sf Mod}$ of the constant sheaf $\pm\breve{\HP}$
over $\text{\sf Ult}(\mathcal{H})$ with stalk $\pm\HP\subset\pm\bdiam\HP_{\mathfrak{u}}$ for all $\mathfrak{u}\in\text{\sf Ult}(\mathcal{H})$.
The image of $\pm\breve{\HP}^{\rm cl}$ (=  the restriction of $\pm\breve{\HP}$ to $\text{\sf Cone}(\mathcal{H})$) in $\bdiam{\sf Mod}$ is denoted 
$\bdiam{\sf Mod}^{\rm cl}$; likewise, the image of the restriction to the union of the $\uptheta$-quantum cone ultrafilters is denoted $\bdiam{\sf Mod}^{\rm qt}$.

Given a hyperfinite set $[F_{\upalpha}]\in\mathcal{H}$ 
and  $\bast z_{[F_{\upalpha}]} = \bast\{ z_{F_{\upalpha}}\}\in\pm\bast\HP\subset\bast\C$, let 
\[ G_{k}(\bast z_{[F_{\upalpha}]}):= \bast\left\{ \sum_{(m,n)\in F_{\upalpha}} (mz_{F_{\upalpha}} +n)^{-2k}\right\}.\]
Following the usual procedure we may then define
$j(\bast z_{[F_{\upalpha}]})$.  Extending to a map of nets indexed by $\mathcal{H}$ leads to a function of sheaves
\[\bdiam j: \pm\bdiam\breve{\HP}\longrightarrow \bdiam\breve{\C} .\]

In the Theorem which follows we will need to specify topologies on the various sheaves and ultrasolenoids defined above.
This can be done by putting a topology on $\bdiam\breve{\C}$ as follows.  First, we topologize each fiber $\bdiam\C_{\mathfrak{u}}$ 
using the $(\bdiam\R_{\mathfrak{u}})_{+}$-valued absolute value $\bdiam |\cdot |_{\mathfrak{u}}$.  Note that there is a canonical
inclusion of the set of $\mathcal{H}$-nets 
\[  \bast \C^{\mathcal{H}}\hookrightarrow \bdiam\breve{\Upgamma} \]
given by 
\[ \{ \bast z_{[F_{\upalpha}]}\}\mapsto  (\mathfrak{u}\mapsto \bdiam \{\bast z_{[F_{\upalpha}]}\}_{\mathfrak{u}}). \]
We denote this section simply $\bdiam z$, and its value at $\mathfrak{u}$ by $\bdiam z_{\mathfrak{u}}$. Now given 
$ \mathcal{O}\subset {\sf Ult}(\mathcal{H})$,  $\{ \bast z_{[F_{\upalpha}]}\}\subset\bast\C$ and
$\{ \bast r_{[F_{\upalpha}]}\}\subset \bast\R_{+}$ a pair of $\mathcal{H}$-nets, let 
\[  \breve{\mathcal{O}}(\bdiam z;\bdiam r) =  \big\{  \bdiam w_{\mathfrak{u}}\big|\; \mathfrak{u}\in \mathcal{O}, | \bdiam z_{\mathfrak{u}}- \bdiam w_{\mathfrak{u}}|_{\mathfrak{u}}<
\bdiam r_{\mathfrak{u}}\big\} .\]
The sets $ \breve{\mathcal{O}}(\bdiam z;\bdiam r)$ form the base for a topology on $\bdiam\breve{\C}$, called
the {\bf {\small ultrasheaf topology}}.  Observe that the subspace topology on $\C\times {\sf Ult}(\mathcal{H})\subset\bdiam\breve{\C}$ coincides with the
product topology.

Note that any section $\bdiam w$ defined by a net 
$\{ \bast w_{[F_{\upalpha}]}\}$ is continuous with respect to the ultrasheaf topology: indeed, if
$\breve{\mathcal{O}}(\bdiam z;\bdiam r)$ contains the point $\bdiam w_{\mathfrak{u}}$, then there exists a subset $X\in\mathfrak{u}$
such that the subnet $\{ \bast w_{[F_{\upalpha}]}\}|_{X}$ is contained in the subnet of balls
$\big\{ B_{\bast r_{[F_{\upalpha}]}}(\bast z_{[F_{\upalpha}]})  \big\}\big|_{X}$, where $B_{\bast r_{[F_{\upalpha}]}}(\bast z_{[F_{\upalpha}]})$
is the ball of radius $\bast r_{[F_{\upalpha}]}$ about $\bast z_{[F_{\upalpha}]}$.  But this implies that this is true for any
ultrafilter $\mathfrak{u'}$ containing $X$.  In other words, the pre-image of the open $\breve{\mathcal{O}}(\bdiam z;\bdiam r)$
is the union of the Stone opens $\mathcal{O}_{X}= \{\mathfrak{u}\ni X\}$, where $X$ is as above.  

More generally,
any map $\bdiam\breve{\C}\longrightarrow \bdiam\breve{\C}$ which takes stalks to stalks is continuous if it is continuous
along the base and if each map $\bdiam\C_{\mathfrak{u}}\rightarrow \bdiam\C_{\mathfrak{u}'}$ is continuous in the norm topologies
In particular, both the shift and the diagonal actions of ${\rm GL}_{2}(\Z )$ on $\pm\bdiam\breve{\HP}\subset\bdiam\breve{\C}$
act by homeomorphisms, properly discontinuously, so that the induced quotient topologies on $\bdiam\hat{\C}$ -- as well as on $\bdiam\widehat{\sf Mod}$
and each of its subsolenoids -- are Hausdorff.  Note that on $\bdiam{\sf Mod}\subset \bdiam\widehat{\sf Mod}$, the subspace topology coincides with
the topology induced by the product topology $\HP\times {\sf Ult}(\mathcal{H})$.


 \begin{theo}\label{solvaluedinvariant}  The map of nets
$   \{  \bast z_{[F_{\upalpha}]} \}_{[F_{\upalpha}]\in\mathcal{H}} \longmapsto \left\{ j(\bast z_{[F_{\upalpha}]})\right\}_{[F_{\upalpha}]\in\mathcal{H}} $
 induces a {\rm continuous} leaf-preserving map of ultrasolenoids
 \[ \bdiam j :  \bdiam \widehat{\sf Mod}\longrightarrow \bdiam\hat{\C}.\]
 The restriction of $\bdiam j$ to $ \bdiam \text{\sf Mod} $ is induced by the universal modular invariant $ \bdiam\hat{\jmath}:\pm\HP\rightarrow \bdiam\hat{\Upgamma}$ by the formula
 \[\bdiam j(z,\mathfrak{u}) =  \bdiam\hat{\jmath}(z)(\mathfrak{u}). \]
 In particular, if we denote by $\bdiam j^{\rm cl}$, $\bdiam j^{\rm qt}$ the restriction of $\bdiam j$ to $\bdiam{\sf Mod}^{\rm cl}$ resp. 
 $\bdiam{\sf Mod}^{\rm qt}$ then 
 \[ \bdiam j^{\rm cl}(z,\mathfrak{u}) =  \bdiam\hat{\jmath}^{\rm cl}(z)(\mathfrak{u}), \quad 
 \bdiam j^{\rm qt}(z,\mathfrak{u}) =  \bdiam\hat{\jmath}^{\rm qt}(z)(\mathfrak{u}).
  \]
 \end{theo}
 
 \begin{proof}  Given $A\in {\rm GL}_{2}(\Z)$, define 
$ A'\circledast \bdiam z_{\mathfrak{u}}\in \bdiam\C_{A^{-T}\mathfrak{u}}$ 
via the action on representative nets:
\[ \left\{A'(\bast z)_{[A^{T}F_{\upalpha}]}\right\}.\]
 Then the map of sheaves  $\bdiam G_{k}: \pm\bdiam\breve{\HP}\longrightarrow \bdiam\breve{\C}$
 satisfies
 \[  ( A'\circledast \bdiam z_{\mathfrak{u}})^{k}\cdot \bdiam G_{k}(A\circledast \bdiam z_{\mathfrak{u}}) = A\odot \bdiam G_{k} ( \bdiam z_{\mathfrak{u}}). \]
 From this, it follows immediately that $\bdiam j$ induces a function
$ \bdiam j :  \bdiam \widehat{\sf Mod}\longrightarrow \bdiam\hat{\C}$, and that the functions $\bdiam j, \bdiam j^{\rm cl}$, $\bdiam j^{\rm qt}$ are related
to $\bdiam \hat{\jmath}$, $ \bdiam\hat{\jmath}^{\rm cl}$, $ \bdiam\hat{\jmath}^{\rm qt}$ by the formulas given in the statement of the Theorem.
 If we fix a net $ \{  \bast z_{[F_{\upalpha}]} \}$, the association $\mathfrak{u}\mapsto \bdiam z_{\mathfrak{u}}$ gives a continuous section of
$\pm\bdiam\breve{\HP}$ and therefore $\mathfrak{u}\mapsto \bdiam j (\bdiam z_{\mathfrak{u}})$ is continuous as well.  If we fix $\mathfrak{u}\in {\sf Ult}(\mathcal{H})$ then the map $\pm\bdiam\HP_{\mathfrak{u}}\rightarrow \bdiam\C_{\mathfrak{u}}$, $\bdiam z_{\mathfrak{u}}\mapsto  \bdiam j (\bdiam z_{\mathfrak{u}})$ gives a continuous map of the stalk $\bdiam\C_{\mathfrak{u}}$: by transference from the continuity of the classical modular invariant.
 Thus $\bdiam j$ is continuous in the ultrasheaf topology.
 \end{proof}
 

Note the dual nature of the classical and quantum invariant: the values of the classical invariant are recovered {\it along
the leaves} of a restriction, whereas the quantum invariant is recovered {\it along the transversal} of another.

\section{Appendix: Experimental Values of $J^{\rm qt}$}

In this section we present some preliminary experimental evidence which suggests that for $\uptheta\in\R$ a quadratic irrationality, $j^{\rm qt}(\uptheta)$
is a finite set.  

Let $D$ be a fundamental discriminant. For $u$ the fundamental unit in $\Q(\sqrt{D})$ the function
\[ x\mapsto J^{\rm qt}_{u^{-x}}(u) \]
is mapped using the {\sf ploth} function of PARI/GP.  
Below are the results for the first five fundamental discriminants $D=5,8,12, 13$ and $17$.

In the case of discrimant $5$, the experimental lower bound $0.8188...$ matches the value obtained by calculating
the formula (\ref{rogramexpress}) and gives $j_{\sf {\small best}}^{\rm qt}(\upvarphi)\approx  9538.2496...$.  On the other hand, if we consider the modification of (\ref{rogramexpress}) obtained
by replacing $G_{M}(\upvarphi)$ by $G_{M}'(\upvarphi)= G_{M}(\upvarphi)+1$, 
a value
very close to $0.8501...$ = the experimental supremum of $J^{\rm qt}(\upvarphi )$ is returned.  

More generally, close examination of the graphs below suggests that the number of returned values of 
$J^{\rm qt}(u)$ is approximately $D$.
Based on these (very preliminary)
computations
 it doesn't seem entirely unreasonable to contemplate the following rough

\begin{conj} Let $\uptheta\in\R-\Q$ be a quadratic irrationality belonging to $\Q (\sqrt{D})$.  Then 
$ j^{\rm qt}(\uptheta )$ is a finite bounded set.  If $\uptheta=u$ is a fundamental unit then
\[ {\rm card}\big(j^{\rm qt}(u )\big) =O(D).\]
\end{conj}

\begin{figure}[htbp]
\centering
\includegraphics[width=5in]{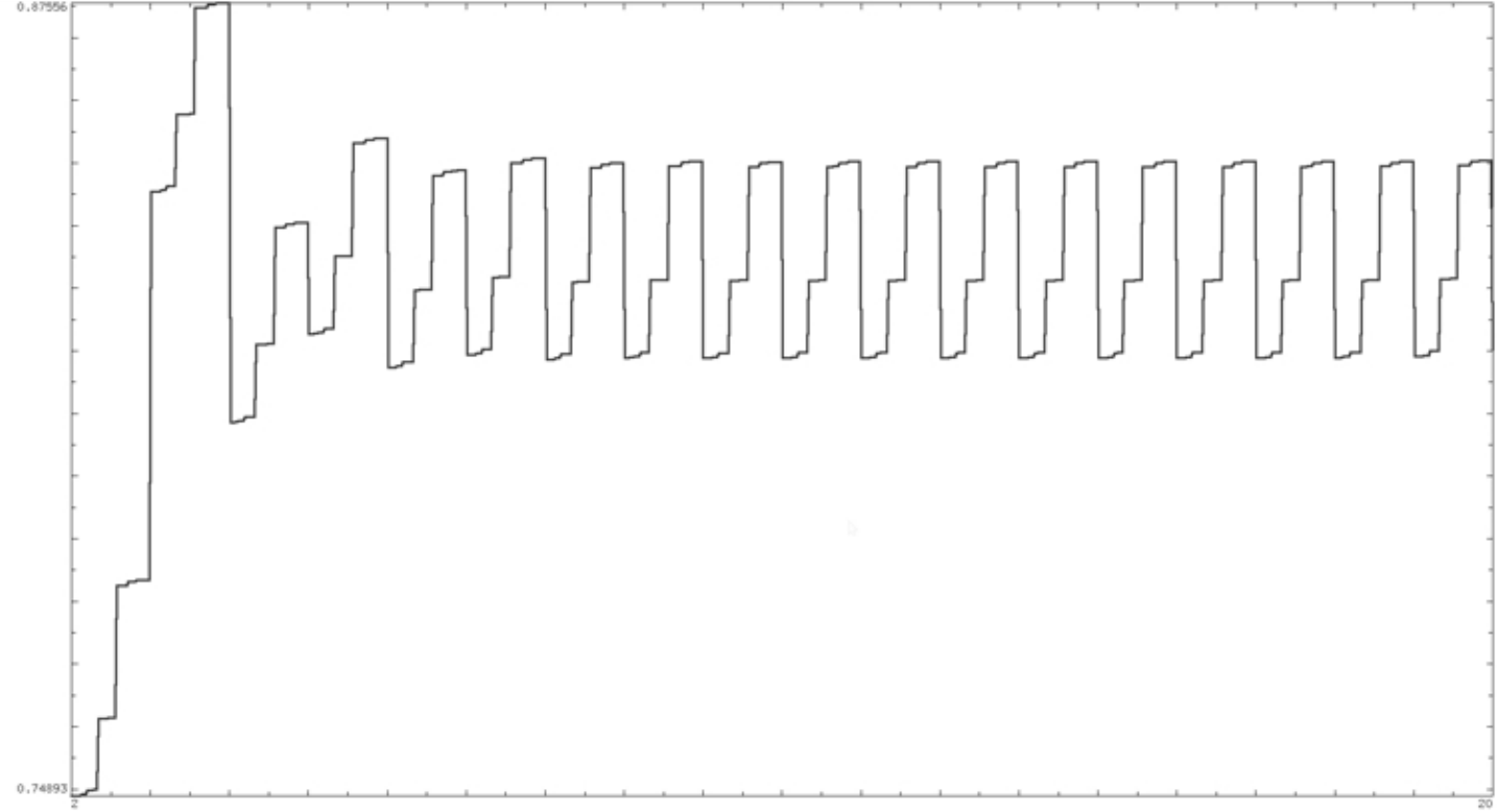}
\caption{D=5.}\label{portraits}
\end{figure}

\begin{figure}[htbp]
\centering
\includegraphics[width=5in]{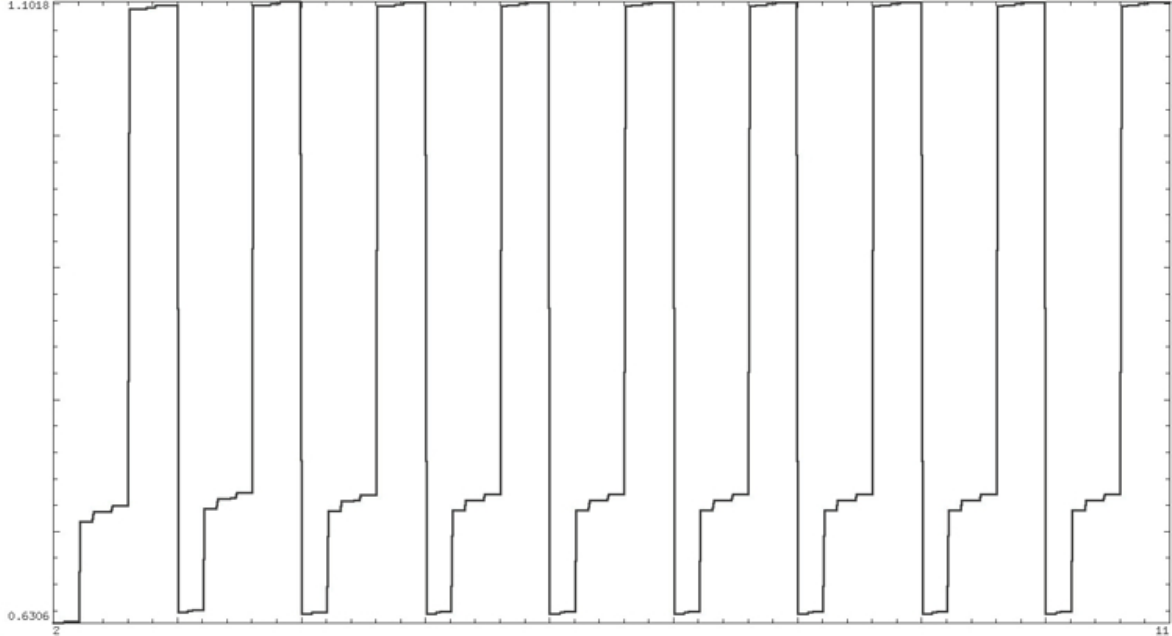}
\caption{D=8.}\label{portraits}
\end{figure}

\begin{figure}[htbp]
\centering
\includegraphics[width=5in]{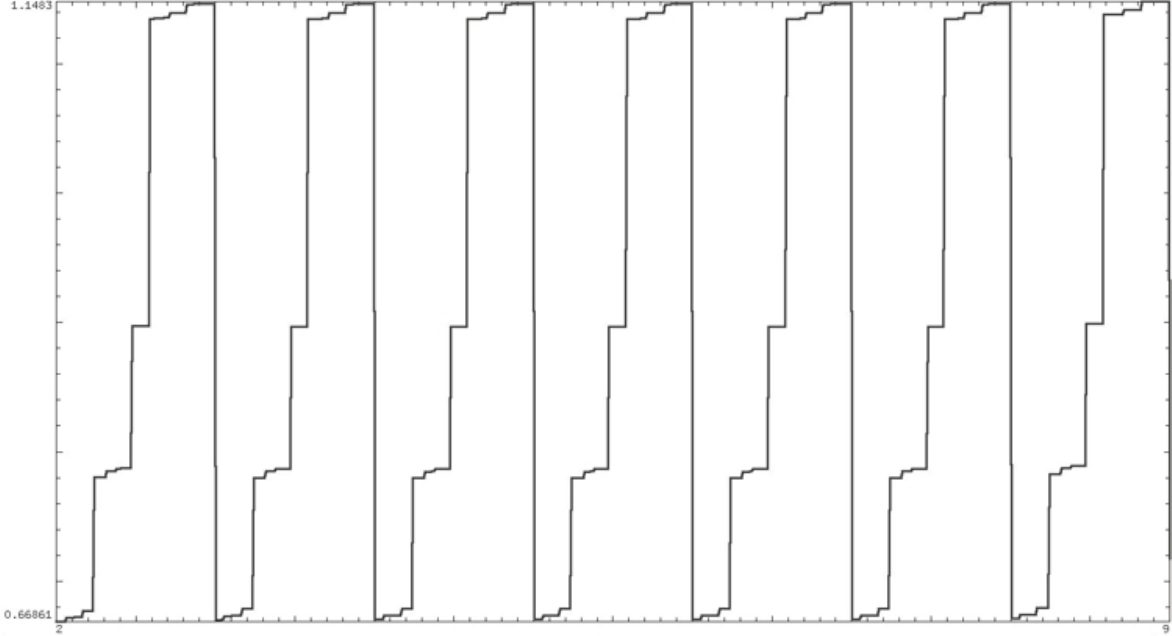}
\caption{D=12.}\label{portraits}
\end{figure}

\begin{figure}[htbp, H]
\centering
\includegraphics[width=5in]{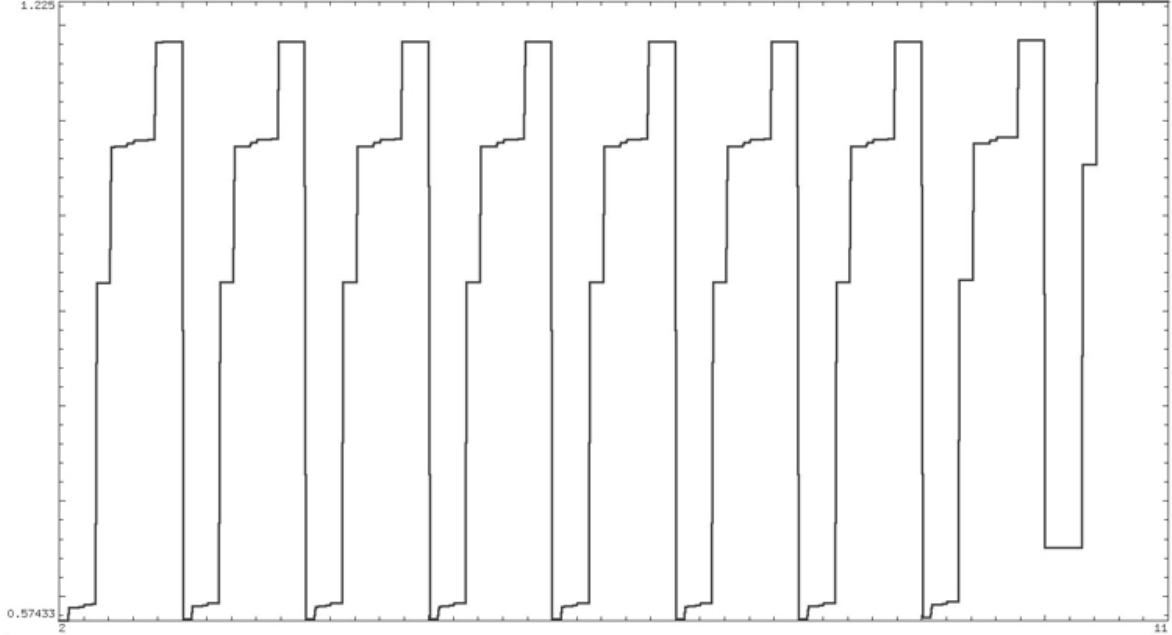}
\caption{D=13.}\label{portraits}
\end{figure}

\newpage

\begin{figure}[ H]
\centering
\includegraphics[width=5in]{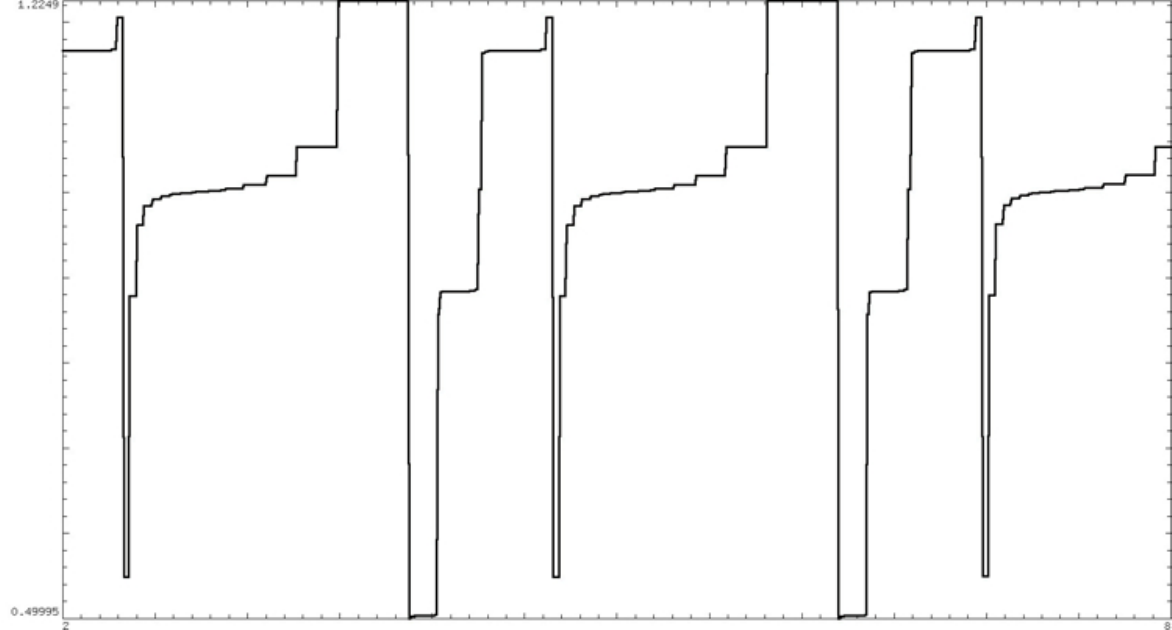}
\caption{D=17.}\label{portraits}
\end{figure}

\newpage

\bibliographystyle{amsalpha}

\end{document}